\documentclass{amsart}
\usepackage{amsmath,amssymb,amsthm,texdraw,citesort}

\numberwithin{equation}{section}

\theoremstyle{plain}
\newtheorem{theorem}{Theorem}[section]
\newtheorem{prop}[theorem]{Proposition}
\newtheorem{lemma}[theorem]{Lemma}
\newtheorem{cor}[theorem]{Corollary}
\theoremstyle{definition}
\newtheorem{definition}[theorem]{Definition}
\newtheorem{example}[theorem]{Example}
\theoremstyle{remark}
\newtheorem{remark}[theorem]{Remark}

\newcommand{\nc}{\newcommand}

\nc{\eit}{\tilde{e}_i}
\nc{\fit}{\tilde{f}_i}
\nc{\ali}{\alpha_i}
\nc{\op}{\oplus}
\nc{\ot}{\otimes}
\nc{\N}{\mathbf{N}}
\nc{\Z}{\mathbf{Z}}
\nc{\Q}{\mathbf{Q}}
\nc{\g}{\mathfrak{g}}
\nc{\uq}{U_q}
\nc{\ft}{\tilde{f}}
\nc{\et}{\tilde{e}}
\nc{\an}{A_n}
\nc{\anone}{A_n^{(1)}}
\nc{\B}{\mathcal{B}}
\nc{\sln}{\mathfrak{sl}_{n+1}}
\nc{\slnh}{\hat{\mathfrak{sl}}_{n+1}}
\nc{\defi}[1]{\emph{\textbf{#1}}}
\nc{\veps}{\varepsilon}
\nc{\vphi}{\varphi}
\nc{\vepst}{\tilde{\varepsilon}}
\nc{\vphit}{\tilde{\varphi}}
\nc{\La}{\Lambda}
\nc{\la}{\lambda}
\nc{\T}{\mathcal{T}}
\nc{\Y}{\mathcal{Y}}
\nc{\spaceiso}{\Psi}
\nc{\spaceosi}{\Phi}
\nc{\M}{\mathcal{M}}
\nc{\Me}{\mathcal{M^{E}}}
\nc{\Minfty}{\mathcal{M^{\infty}}}
\nc{\hwc}{\mathcal{B}}
\nc{\Binf}{\mathcal{B}(\infty)}
\nc{\Tinf}{\mathcal{T}(\infty)}

\DeclareMathOperator{\wt}{wt}
\DeclareMathOperator{\wtt}{\tilde{wt}}

\begin{document}

\title
[Descriptions of the crystal $\B(\infty)$ for $G_2$]
{Descriptions of the crystal $\B(\infty)$ for $G_2$}

\author{Hyeonmi Lee}
\address{Korea Institute for Advanced Study\\
         207-43 Cheongnyangni 2-dong, Dongdaemun-gu\\
         Seoul 130-722, Korea}
\email{hmlee@kias.re.kr}
\thanks{This work was supported in part by KOSEF Grant
        R01-2003-000-10012-0}

\begin{abstract}
We study the crystal base of the negative part of a quantum group.
Two explicit descriptions of the crystal $\B(\infty)$ for types
$G_2$ are given. The first is given in terms of \emph{extended
Nakajima monomials} and the second realization follows a similar
result given for other finite types by \emph{Cliff}.
\end{abstract}

\maketitle

%%%%%%%%%%%%%%%%%%%%%%%%%%%%%%%%%%%%%%%%%%%%%%%%%%%%%%%%%%%%%%%%%%%%%%%%%%
\section{Introduction}%
\label{Introduction}

Quantum group $U_q(\mathfrak{g})$ is a $q$-deformation of the
universal enveloping algebra over a Lie algebra $\mathfrak{g}$,
and crystal bases reveal the structure of
$U_q(\mathfrak{g})$-modules in a very simplified form. As these
$U_q(\mathfrak{g})$-modules are known to be $q$-deformations of
modules over the original Lie algebras, knowledge of these
structures also affects the study of Lie algebras.

The crystal $\Binf$, which is the crystal base of the negative
part $U_q^-(\mathfrak{g})$ of a quantum group, has received
attention since the very birth of crystal base
theory~\cite{K91,K90}. This is not only because it is an essential
part of the \emph{grand loop} argument proving the existence of
crystal bases, but because it gives insight into the structure of
quantum group itself.

Much effort has been
made~\cite{MR1475048,ks,saito,MR1614241,lee,lee2,HL} to give
explicit description of the crystals $\Binf$ over various
Kac-Moody algebras. A related known result is that it is possible
to characterize the highest weight crystal $\B(\la)$ over
symmetrizable Kac-Moody algebras, in terms of Nakajima
monomials~\cite{MR1988989,MR2043368,MR2074700,kim,shin}, an object
which was introduced by Nakajima~\cite{MR1865400,MR1988990}. This
has lead to the belief that it should be possible to give a
similar description for $\Binf$ also. Starting from a theorem of
Kashiwara and Nakajima on the crystal structure of
monomials~\cite{MR1988989}, we can argue that it is not possible
to find the crystal $\B(\infty)$ within the set of Nakajima
monomials with their given crystal structure. Hence, in our recent
work~\cite{lee}, we constructed the set of extended monomials and
developed a crystal structure on it, conjecturing that a certain
connected component of the crystal would be isomorphic to
$\B(\infty)$. Actually, the set of Nakajima monomials can be
embedded as a subcrystal in this set of extended Nakajima
monomials. Thus, the monomial theory developed for irreducible
highest weight crystal can easily be transferred to the extended
monomial set.

In the current work, we restrict ourselves to the $G_2$-type
finite simple Lie algebra. For this case, we give an explicit set
of extended Nakajima monomials and show it to be isomorphic to
$\Binf$. The previous work~\cite{HL}, giving a Young tableaux
realization of $\Binf$ for the finite simple types, is used in
doing this. We also extend Cliff's~\cite{MR1614241} realization of
$\Binf$ for classical finite types, given in terms of a completely
different object, to the $G_2$-type.

The paper is organized as follows. We start by reviewing the
notion of extended Nakajima monomials and the crystal structure
given on the set of such monomials. Also, we cite Young tableau
expression of crystal $\B(\infty)$ for type $G_2$ which play a
crucial role in our work. We then proceed to give a monomial
realization of the crystal $\B(\infty)$. In the process of
obtaining these results, we give new expressions for the Kashiwara
operators acting on the extended Nakajima monomials, more
appropriate for the situation in hand. In the last section, we
deal with Cliff's approach of realizing $\Binf$.

%%%%%%%%%%%%%%%%%%%%%%%%%%%%%%%%%%%%%%%%%%%%%%%%%%%%%%%%%%%%%%%%%%%%%%%%%%
\section{Extended Nakajima monomials and Young tableaux}%
\label{Nakajima's monomials and crystals}

In this section, we introduce notation and cite facts that are crucial
for our work. Please refer to the references cited in the introduction
or books on quantum groups~\cite{MR1881971,MR1359532} for the basic
concepts on quantum groups and crystal bases.

Let us first fix the basic notation.
\begin{itemize}
\item $I=\{1,2\}$ : index set for $G_2$-type. \item
$A=(a_{ij})_{i,j\in I}$ : Cartan matrix of type $G_2$. \item
$\ali$, $\La_i$ $(i\in I)$ : simple root, fundamental weight.
\item $\Pi =\{\ali\vert i\in I\}$ : the set of simple root. \item
$P =\op_{i\in I} \Z\La_i$ : weight lattice. \item $\uq(G_2)$ :
quantum group for $G_2$. \item $\uq^-(G_2)$ : subalgebra of
$\uq(G_2)$
      generated by $f_i$ ($i\in I$).
\item $\fit, \eit$ : Kashiwara operators. \item $\hwc(\infty)$ :
crystal base of $\uq^-(G_2)$.
\end{itemize}
Throughout this paper, a $U_q(G_2)$-crystal will refer to a
(abstract) crystal associated with the Cartan datum $(A,\Pi,P)$.

%%%%%%%%%%%%%%%%%%%%%%%%%%%%%%%%%%%%%%%%%%%%
\subsection{Nakajima monomials}%
\label{Nakajima monomials}

We now recall the set of monomials and its crystal structure
discovered by Nakajima~\cite{MR1988990} and also recall their
extension introduced in~\cite{lee}. Both of these sets were
defined for all symmetrizable Kac-Moody algebras, but we shall
restrict ourselves to the $G_2$ case in this paper.

Let $\Me$ be a certain set of formal monomials in the variables
$Y_i(m)$ ($i\in I$, $m\in\Z$). More explicitly,
\begin{equation}
\Me =
\left\{\prod_{(i,m)\in I\times\Z}
        {Y_i(m)}^{y_i(m)}
        \Big\vert
\begin{aligned}
 & y_i(m)\!=(y^0_i(m),y^1_i(m)) \in\Z\times\Z \, \textup{vanishes}\\
 & \textup{except at finitely many} \ (i,m)
\end{aligned}
\right\}.
\end{equation}
We give the lexicographic order to the set $\Z\times\Z$ of variable
exponents. Fix any set of integers $c={(c_{ij})}_{i \neq j \in I}$ such
that
\begin{equation}
c_{ij}+c_{ji}=1,
\end{equation}
and set
\begin{equation}
A_i(m)={Y_i(m)}^{(0,1)}{Y_i(m+1)}^{(0,1)}
\prod_{j \neq i} {Y_j(m+c_{ji})}^{(0,\langle h_j,\alpha_i\rangle)}.
\end{equation}

The crystal structure on $\Me$ is defined as follows. For every
monomial $M=\prod_{(i,m)\in I\times\Z}{Y_i(m)}^{y_i(m)}$, we set
\begin{align}
\wtt(M) &= \sum_i(\sum_m y_i(m))\La_i,\\
\vphit_i(M) &= \textup{max}\{ \sum_{k\le m} y_i(k)
            \,\Big\vert \ m\in\Z \},\\
\vepst_i(M) &= \textup{max}\{ -\sum_{k>m} y_i(k)
            \,\Big\vert \ m\in\Z \}.
\end{align}
Notice that the coefficients of $\wtt(M)$ are pairs of integers. In
this setting, we have $\vphit_i(M)\geq (0,0)$, $\vepst_i(M)\geq (0,0)$,
and $\wtt(M)=\sum_i (\vphit_i(M)-\vepst_i(M))\La_i$. Set
\begin{alignat}{2}
\wt(M) &= \sum_i(\sum_m y^1_i(m))\La_i,\label{structure1}\\
\vphi_i(M) &= \sum_{k\le m} y^1_i(k) & &\textup{where}\
              \vphit_i(M)= \sum_{k\le m} (y^0_i(k),y^1_i(k)),
              \label{structure2}\\
\veps_i(M) &= -\sum_{k>m} y^1_i(k) & &\textup{where}\
              \vepst_i(M)= -\sum_{k>m} (y^0_i(k),y^1_i(k)).
              \label{structure3}
\end{alignat}
Then we trivially have $\wt(M)=\sum_i (\vphi_i(M)-\veps_i(M))\La_i$.
 From the above definition, ${Y_i(m)}^{(0,1)}$ has the weight $\La_i$,
and so $A_i(m)$ has the weight $\alpha_i$. We define the action of
Kashiwara operators by
\begin{align}
\fit(M) &=
\begin{cases}\label{actions}
0 &\textup{if} \ \vphit_i(M)=(0,0),\\
A_i(m_f)^{-1}M &\textup{if} \ \vphit_i(M)>(0,0),
\end{cases}\\
\eit(M) &=
\begin{cases}\label{actions2}
0 &\textup{if} \ \vepst_i(M)=(0,0),\\
A_i(m_e)M &\textup{if} \ \vepst_i(M)>(0,0).
\end{cases}
\end{align}
Here,
\begin{align}
m_f &=\textup{min}\{m | \vphit_i(M)=
      \sum_{k\le m} y_i(k) \}
    =\textup{min}\{m | \vepst_i(M)=
      -\sum_{k>m} y_i(k) \},\label{emf}\\
m_e &=\textup{max}\{m | \vphit_i(M)=
      \sum_{k\le m} y_i(k) \}
    =\textup{max}\{m | \vepst_i(M)=
      -\sum_{k>m} y_i(k) \}.\label{eme}
\end{align}
Note that $y_i(m_f)>(0,0)$, $y_i(m_f+1)\le (0,0)$,
$y_i(m_e+1)<(0,0)$, and $y_i(m_e)\ge (0,0)$.

The Kashiwara operators, together with the maps
$\vphi_i$, $\veps_i$ $(i\in I)$, $\wt$, define
a crystal structure on the set $\Me$~\cite{lee}.

The set of monomials $\prod_{(i,m)\in
I\times\Z}{Y_i(m)}^{(y^0_i(m),y^1_i(m))}$ of $\Me$ with $y^0_i(m)=0$ is
exactly the Nakajima monomial set $\M$ if we identify
${Y_i(m)}^{(0,y^1_i(m))}\in\Me$ with ${Y_i(m)}^{y^1_i(m)}\in\M$. The
crystal structure on $\M$, introduced in~\cite{MR1988989}, is
compatible with that on $\Me$ under this identification. The crystal
$\M$ is a subcrystal of $\Me$.

Restriction of the following theorem to just the monomials of $\M$
appears in~\cite{MR1988989}.

\begin{theorem}\textup{(}\cite{lee}\textup{)}\label{mla}
If a monomial $M\in\Me$ of $\wtt(M)=\sum_i(0,p_i)\La_i$, where each
$p_i$ is a nonnegative integer, satisfies $\eit(M)=0$ for all $i\in I$,
then the connected component of $\Me$ containing $M$ is isomorphic to
$\B(\sum_i p_i\La_i)$ as a $\uq(\g)$-crystal.
Conversely, given any subset of $\Me$ isomorphic to $\B(\sum_i p_i\La_i)$,
there exists an element $M$ in the subset such that
$\wtt(M)=\sum_i(0,p_i)\La_i$ and $\eit(M)=0$ for all $i\in I$.
\end{theorem}

The following is a conjecture on the crystal $\B(\infty)$ introduced
in~\cite{lee} and stated for all symmetrizable Kac-Moody algebras. Its
converse is known to be true~\cite{lee}.

\vspace{2mm}
\emph{If a monomial $M\in\Me$ of $\wtt(M)=\sum_i(p_i,0)\La_i$, where each
$p_i$ is a positive integer, satisfies $\eit(M)=0$ for all $i\in I$,
then the connected component of $\Me$ containing $M$ is isomorphic to
$\B(\infty)$ as a $\uq(\g)$-crystal.}
\vspace{2mm}

In~\cite{lee,lee2}, result for $\B(\infty)$ of type $\an$ and
$\anone$ was given as evidence supporting this conjecture. In the
next section, we will give a concrete listing of elements
containing a weight zero vector $M\in\Me$ mentioned in the above
theorem, for the case of $G_2$. We shall show this set to be a
crystal and give an isomorphism between the crystal and crystal
$\B(\infty)$. This is another result supporting the above
conjecture.

%%%%%%%%%%%%%%%%%%%%%%%%%%%%%%%%%%%%%%%%%%%%
\subsection{Young tableaux}%
\label{Young tableaux}

In this section, we recall a Young tableaux description for the crystal
$\B(\infty)$ over type $G_2$ introduced in~\cite{HL}.

For the $G_2$-type, we shall take the Young tableau realization of
highest weight crystal $\B(\la)$ given in~\cite{KM} as the definition
of semi-standard tableaux. Since the work is a rather well known
result, we refer readers to the original papers and shall not repeat
the complicated definition here. The alphabet to be used inside the
boxes constituting the Young tableaux will be denoted by $J$, and it
will be equipped with an ordering $\prec$, as given in~\cite{KM}.

\begin{equation*}
J=\{1\prec{2}\prec{3}\prec{0}\prec\bar{3}\prec\bar{2}\prec\bar{1}\}.
\end{equation*}

\begin{definition}\hfill
\begin{enumerate}
\item A semi-standard tableau $T$ of shape $\la\in P^+$,
      equivalently, an element of an irreducible highest weight
      crystal $\B(\la)$ for the $G_2$ type, is \emph{large} if it
      consists of $2$ non-empty rows, and if the
      number of $1$-boxes in the first row is strictly greater
      than the number of all boxes in the second row and
      the second row contains at least one $2$-box.
\item A large tableau $T$ is \emph{marginally} large,
      if the number of $1$-boxes in the first row of $T$ is greater than
      the number of all boxes in the second row by exactly one and
      the second row of $T$ contain one $2$-box.
\end{enumerate}
\end{definition}

In Figure~\ref{tbl:1}, we give examples of semi-standard tableaux.
The one on the left is large, the one on the middle is marginally
large, and the one on the right is not large.
\begin{figure}
\centering
\begin{tabular}{rll}
 \raisebox{-0.4\height}{
 \begin{texdraw}%
 \drawdim in
 \arrowheadsize l:0.065 w:0.03
 \arrowheadtype t:F
 \fontsize{6}{6}\selectfont
 \textref h:C v:C
 \drawdim em
 \setunitscale 1.4
 \move(0 2)\rlvec(10 0)
 \move(0 1)\rlvec(10 0)\rlvec(0 1)
 \move(0 0)\rlvec(3 0)\rlvec(0 2)
 \move(0 0)\rlvec(0 2)
 \move(1 0)\rlvec(0 2)
 \move(2 0)\rlvec(0 2)
 \move(6 1)\rlvec(0 1)
 \move(7 1)\rlvec(0 1)
 \move(8 1)\rlvec(0 1)
 \move(9 1)\rlvec(0 1)
 \move(5 1)\rlvec(0 1)
 \move(4 1)\rlvec(0 1)
 \htext(0.5 1.5){$1$}
 \htext(1.5 1.5){$1$}
 \htext(2.5 1.5){$1$}
 \htext(3.5 1.5){$1$}
 \htext(9.5 1.5){$\bar 1$}
 \htext(8.5 1.5){$\bar 2$}
 \htext(7.5 1.5){$\bar 3$}
 \htext(6.5 1.5){$0$}
 \htext(5.5 1.5){$3$}
 \htext(4.5 1.5){$2$}
 \htext(0.5 0.5){$2$}
 \htext(1.5 0.5){$2$}
 \htext(2.5 0.5){$3$}
 \htext(1.5 0.5){$2$}
 \htext(2.5 0.5){$3$}
 \end{texdraw}%
 }
 & \raisebox{-0.4\height}{
 \begin{texdraw}%
 \drawdim in
 \arrowheadsize l:0.065 w:0.03
 \arrowheadtype t:F
 \fontsize{6}{6}\selectfont
 \textref h:C v:C
 \drawdim em
 \setunitscale 1.4
 \move(1 2)\rlvec(9 0)
 \move(1 1)\rlvec(9 0)\rlvec(0 1)
 \move(1 0)\rlvec(2 0)\rlvec(0 2)
 \move(1 0)\rlvec(0 2)
 \move(2 0)\rlvec(0 2)
 \move(6 1)\rlvec(0 1)
 \move(7 1)\rlvec(0 1)
 \move(8 1)\rlvec(0 1)
 \move(9 1)\rlvec(0 1)
 \move(5 1)\rlvec(0 1)
 \move(4 1)\rlvec(0 1)
 \htext(1.5 1.5){$1$}
 \htext(2.5 1.5){$1$}
 \htext(3.5 1.5){$1$}
 \htext(9.5 1.5){$\bar 1$}
 \htext(8.5 1.5){$\bar 2$}
 \htext(7.5 1.5){$\bar 3$}
 \htext(6.5 1.5){$0$}
 \htext(5.5 1.5){$3$}
 \htext(4.5 1.5){$2$}
 \htext(1.5 0.5){$2$}
 \htext(2.5 0.5){$3$}
 \htext(1.5 0.5){$2$}
 \htext(2.5 0.5){$3$}
 \end{texdraw}%
 }
 &
 \raisebox{-0.4\height}{
 \begin{texdraw}%
 \drawdim in
 \arrowheadsize l:0.065 w:0.03
 \arrowheadtype t:F
 \fontsize{6}{6}\selectfont
 \textref h:C v:C
 \drawdim em
 \setunitscale 1.4
 \move(0 2)\rlvec(6 0)
 \move(0 1)\rlvec(6 0)\rlvec(0 1)
 \move(0 0)\rlvec(3 0)\rlvec(0 2)
 \move(0 0)\rlvec(0 2)
 \move(1 0)\rlvec(0 2)
 \move(2 0)\rlvec(0 2)
 \move(5 1)\rlvec(0 1)
 \move(4 1)\rlvec(0 1)
 \htext(0.5 1.5){$1$}
 \htext(1.5 1.5){$2$}
 \htext(2.5 1.5){$0$}
 \htext(3.5 1.5){$\bar{3}$}
 \htext(4.5 1.5){$\bar{3}$}
 \htext(5.5 1.5){$\bar{2}$}
 \htext(0.5 0.5){$2$}
 \htext(1.5 0.5){$3$}
 \htext(2.5 0.5){$\bar{2}$}
 \end{texdraw}%
 }
\end{tabular}\\[1em]
\caption{Large (left), marginally large (middle), and non-large
(right) tableaux}\label{tbl:1}
\end{figure}

\begin{definition}
We denote by $\T(\infty)$ the set of all marginally large
tableaux. The marginally large tableau whose $i$-th row consists
only of $i$-boxes ($i\in I$) is denoted by $T_{\infty}$.
\end{definition}

The set $\Tinf$, consists of all tableaux of the following form.
The unshaded part must exist, whereas the shaded part is optional
with variable size.
\begin{equation*}
T = \raisebox{-0.5\height}{\ %
\begin{texdraw}
 \fontsize{6}{6}\selectfont
 \textref h:C v:C
 \drawdim em
 \setunitscale 1.35
 \move(14.5 0)\lvec(1 0)\lvec(1 -1)\lvec(14.5 -1)\lvec(14.5 0)
 \ifill f:0.8
 \move(0 0)\lvec(-2.5 0)\lvec(-2.5 -2)\lvec(0 -2)\lvec(0 0)
 \ifill f:0.8
 \move(14.5 0)\lvec(-3.5 0)
 \move(14.5 -1)\lvec(-3.5 -1)
 \move(0 -2)\lvec(-3.5 -2)
 \move(-3.5 0)\lvec(-3.5 -2)
 \htext(-3 -0.5){$1$}
 \htext(-0.55 -0.5){$1$}
 \htext(-1.7 -0.5){$\cdots$}

 \move(0 0)\rlvec(0 -2)
 \htext(0.5 -0.5){$1$}
 \move(1 0)\rlvec(0 -1)
 \htext(2.25 -0.5){$2\!\cdots\!2$}
 \move(3.5 0)\rlvec(0 -1)
 \htext(4.75 -0.5){$3\!\cdots\!3$}
 \move(6 0)\rlvec(0 -1)
 \htext(6.5 -0.5){$0$}
 \move(7 0)\rlvec(0 -1)
 \htext(8.25 -0.5){$\bar{3}\!\cdots\!\bar{3}$}
 \move(9.5 0)\rlvec(0 -1)
 \htext(10.75 -0.5){$\bar{2}\!\cdots\!\bar{2}$}
 \move(12 0)\rlvec(0 -1)
 \htext(13.25 -0.5){$\bar{1}\!\cdots\!\bar{1}$}
 \move(14.5 0)\rlvec(0 -1)
 \htext(-1.25 -1.5){$3\!\cdots\!3$}
 \move(-2.5 -0.9)\rlvec(0 -1.1)
 \htext(-3 -1.5){$2$}
\end{texdraw}}.
\end{equation*}
The element $T_{\infty}$ is
\begin{equation*}
T_{\infty} = \raisebox{-0.5\height}{\ %
\begin{texdraw}
 \fontsize{6}{6}\selectfont
 \textref h:C v:C
 \drawdim em
 \setunitscale 1.35
 \move(0 0)\lvec(2 0)\lvec(2 -1)\lvec(0 -1)
 \move(1 0)\lvec(1 -2)\lvec(0 -2)\lvec(0 0)
 \htext(0.5 -0.5){$1$}
 \htext(1.5 -0.5){$1$}
 \htext(0.5 -1.5){$2$}
%\drawbb
\end{texdraw}}\,.
\end{equation*}

We recall the action of Kashiwara operators $\fit$, $\eit$ $(i\in I)$
on marginally large tableaux $T\in\T(\infty)$.
\begin{enumerate}
\item We first read the boxes in the tableau $T$ through
      the \emph{far eastern reading} and write down the boxes
      in \emph{tensor product form.}
      That is, we read through each column from top to bottom starting
      from the rightmost column, continuing to the left, and lay down
      the read boxes from left to right in tensor product form.
\item Under each tensor component $x$ of $T$,
      write down $\veps_i(x)$-many 1s followed by $\vphi_i(x)$-many 0s.
      Then, from the long sequence of mixed 0s and 1s,
      successively cancel out every occurrence of (0,1) pair
      until we arrive at a sequence of 1s followed by 0s,
      reading from left to right. This is called the
      $i$-signature of $T$.
\item Denote by $T'$, the tableau obtained from $T$,
      by replacing the box $x$ corresponding to the leftmost $0$
      in the $i$-signature of $T$ with the box $\fit x$.
      \begin{itemize}
      \item If $T'$ is a large tableau,
            it is automatically marginally large.
            We define $\fit T$ to be $T'$.
      \item If $T'$ is not large, then we define $\fit T$ to be the
            large tableau obtained by inserting one column
            consisting of $i$ rows to the left of the box $\fit$
            acted upon. The added column should have a
            $k$-box at the $k$-th row for $1\leq k \leq i$.
      \end{itemize}
\item Denote by $T'$, the tableau obtained from $T$,
      by replacing the box $x$
      corresponding to the rightmost $1$ in the $i$-signature of
      $T$ with the box $\eit x$.
      \begin{itemize}
      \item If $T'$ is a marginally large tableau,
            then we define $\eit T$ to be $T'$.
      \item If $T'$ is large but not marginally large,
            then we define $\eit T$ to be the
            large tableau obtained by removing the column
            containing the changed box.
            It will be of $i$ rows and have a $k$-box
            at the $k$-th row for $1\leq k \leq i$.
      \end{itemize}
\item If there is no $1$ in the $i$-signature of $T$,
      we define $\eit T=0$.
\end{enumerate}

Let $T$ be a tableau in $\T(\infty)$ with the second row
consisting of $b^2_3$-many $3$\,s, one $2$ and the first row
consisting of $b^1_j$-many $j$\,s ($1\prec j\preceq \bar 1$),
$(b^2_3+2)$-many $1$\,s. We define the maps $\wt: \T(\infty)
\rightarrow P$, $\vphi_i, \veps_i : \T(\infty) \rightarrow \Z$ by
setting
\begin{align}
\wt(T)      &=(-b_2^1-b_3^1-2b_0^1-3b_{\bar 3}^1-3b_{\bar 2}^1
              -4b_{\bar 1}^1)\alpha_1 \label{same}\\
            &\quad+(-b_3^1-b_0^1-b_{\bar 3}^1-2b_{\bar 2}^1
              -2b_{\bar 1}^1-b_3^2)\alpha_2,\notag\\
\veps_i(T)  &=\text{the number of $1$s
             in the $i$-signature of $T$},\label{same2}\\
\vphi_i(T)  &=\veps_i(T)+\langle h_i,\wt(T)\rangle.\label{same3}
\end{align}

\begin{theorem}\textup{(}\cite{HL}\textup{)}\label{HL}
The operator given by equations~\eqref{same} to~\eqref{same3}, with
Kashiwara operators define a crystal structure on $\Tinf$. And the
crystal $\T(\infty)$ is isomorphic to $\B(\infty)$ as a
$U_q(G_2)$-crystal.
\end{theorem}

%%%%%%%%%%%%%%%%%%%%%%%%%%%%%%%%%%%%%%%%%%%%%%%%%%%%%%%%%%%%%%%%%%%%%%%%%%
\section{Monomial description of $\B(\infty)$}%
\label{Extended monomial description}

We give a new realization of the crystal $\B(\infty)$, for $G_2$-type,
in terms of extended monomials.

For simplicity, from now on, we take the set $C=(c_{ij})_{i\neq j \in I}$
to be $c_{12}=1$ and $c_{21}=0$.
Then for $m\in\Z$, we have
\begin{equation}
\begin{cases}
A_1(m)={Y_1(m)}^{(0,1)}{Y_1(m+1)}^{(0,1)}{Y_2(m)}^{(0,-1)},\\
A_2(m)={Y_2(m)}^{(0,1)}{Y_2(m+1)}^{(0,1)}{Y_1(m+1)}^{(0,-3)}.
\end{cases}
\end{equation}

The set we define below was originally obtained by applying Kashiwara
actions $\fit$ continuously on the single element
${Y_1(-1)}^{(1,0)}{Y_2(-2)}^{(1,0)}\in\Me$. This choice of starting
monomial will allow us to relate monomials of the set defined below to
tableaux in $\T(\infty)$ naturally.

\begin{definition}
Consider elements of $\Me$ having the form
\begin{equation}\label{extend}
\begin{aligned}
      M=&{Y_1(-1)}^{(1,a_1^{-1})}{Y_1(0)}^{(0,a_1^0)}
        {Y_1(1)}^{(0,a_1^1)}{Y_1(2)}^{(0,a_1^2)}\\
        &\cdot {Y_2(-2)}^{(1,a_2^{-2})}{Y_2(-1)}^{(0,a_2^{-1})}
        {Y_2(0)}^{(0,a_2^0)}{Y_2(1)}^{(0,a_2^1)}
\end{aligned}
\end{equation}
with conditions
\begin{enumerate}
\item $(a_2^{-2}-a_2^{-1})$, $a_2^1$, $a_1^2$, $a_2^{-2}\le 0$,
\item $(a_1^{-1}-a_1^1-a_1^2)+(2a_2^{-2}+a_2^{-1}-a_2^0-2a_2^1)=0$ and\\
      $(a_1^{-1}+a_1^0-a_1^2)+(a_2^{-2}+2a_2^{-1}+a_2^0-a_2^1)=0$,
\item $(a_1^0+a_2^{-1}-a_2^{-2}), (-a_1^1-a_2^1)\in 2\Z_{\ge 0}$ or\\
      $(a_1^0+a_2^{-1}-a_2^{-2}), (-a_1^1-a_2^1)\in \Z_{\ge 0}$ and odd.
\end{enumerate}
Specifically, in case of $a_i^j=0$ for all $i,j$, we have
\begin{equation}\label{hi9}
      M={Y_1(-1)}^{(1,0)}{Y_2(-2)}^{(1,0)}.
\end{equation}
We denote by $\M(\infty)$ the set of all monomials of these form
and by $M_{\infty}$ the monomial of~\eqref{hi9}.
\end{definition}

\begin{figure}
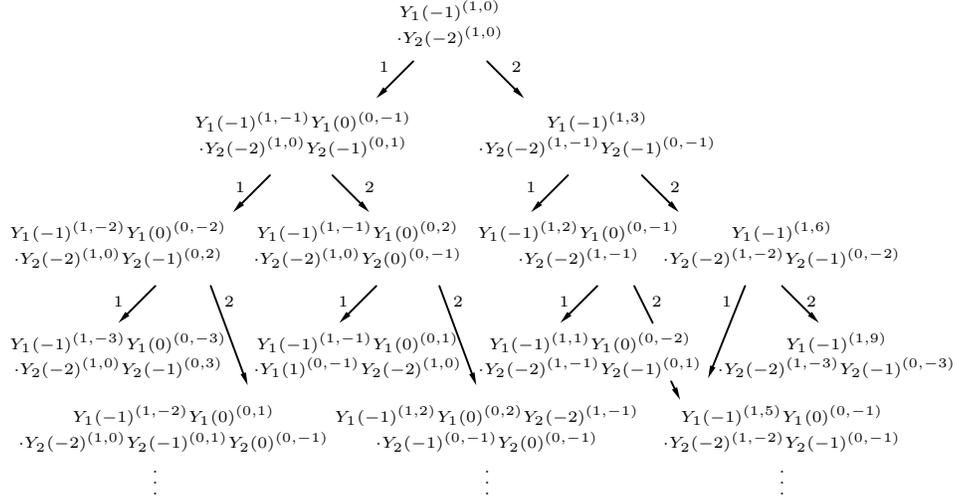

\centering
\begin{texdraw}%
\drawdim in
\arrowheadsize l:0.065 w:0.03
\arrowheadtype t:F
\fontsize{6}{6}\selectfont
\textref h:C v:C
\drawdim em
\setunitscale 1.9
\move(-1 -1.4)\ravec(-1 -1)%%%%%%%%%%%%%%%%%
\move(1 -1.4)\ravec(1 -1)%%%%%%%%%%%%%%%%%
\move(-3.2 -4.5)\ravec(1 -1)%%%%%%%%%%%%%%%%%
\move(-4.9 -4.5)\ravec(-1 -1)%%%%%%%%%%%%%%%%%
\move(3.2 -4.5)\ravec(-1 -1)%%%%%%%%%%%%%%%%%
\move(5.2 -4.5)\ravec(1 -1)%%%%%%%%%%%%%%%%%
\move(-8 -7.5)\ravec(-1 -1)%%%%%%%%%%%%%%%%%
\move(-6.5 -7.5)\ravec(1 -2.7)%%%%%%%%%%%%%%%%%
\move(-2 -7.5)\ravec(-1 -1)%%%%%%%%%%%%%%%%%
\move(-0.3 -7.5)\ravec(1 -2.7)%%%%%%%%%%%%%%%%%
\move(4 -7.5)\ravec(-1 -1)%%%%%%%%%%%%%%%%%
\move(5 -7.5)\rlvec(0.5 -1)%%%%%%%%%%%%%%%%%
\move(6 -10.1)\ravec(0.25 -0.5)%%%%%%%%%%%%%%%%%
\move(8 -7.5)\ravec(-1 -2.5)%%%%%%%%%%%%%%%%%
\move(9 -7.5)\ravec(1 -1)%%%%%%%%%%%%%%%%%
\htext(-8 -12.6){$\vdots$}%%%%%%%%%%%%%%%%%
\htext(9 -12.6){$\vdots$}%%%%%%%%%%%%%%%%%
\htext(1 -12.6){$\vdots$}%%%%%%%%%%%%%%%%%
%%%%%%%%%%%%%%%%%%%%%%%%%%%%%%%%%%%%%%%%%%%%%%%
\htext(0 0){${Y_1(-1)}^{(1,0)}$}
\htext(0 -0.7){$\cdot {Y_2(-2)}^{(1,0)}$}
%%%%%%%%%%%%%%%%%%%%%%%%%%%%%%%%%%%%%%%%%%%%%%%
\htext(-4 -3){${Y_1(-1)}^{(1,-1)}{Y_1(0)}^{(0,-1)}$}
\htext(-4 -3.7){$\cdot {Y_2(-2)}^{(1,0)}{Y_2(-1)}^{(0,1)}$}

\htext(4 -3){${Y_1(-1)}^{(1,3)}$}
\htext(4 -3.7){$\cdot {Y_2(-2)}^{(1,-1)}{Y_2(-1)}^{(0,-1)}$}
%%%%%%%%%%%%%%%%%%%%%%%%%%%%%%%%%%%%%%%%%%%%%%
\htext(-9 -6){${Y_1(-1)}^{(1,-2)}{Y_1(0)}^{(0,-2)}$}
\htext(-9 -6.7){$\cdot {Y_2(-2)}^{(1,0)}{Y_2(-1)}^{(0,2)}$}

\htext(-2.5 -6){${Y_1(-1)}^{(1,-1)}{Y_1(0)}^{(0,2)}$}
\htext(-2.5 -6.7){$\cdot {Y_2(-2)}^{(1,0)}{Y_2(0)}^{(0,-1)}$}

\htext(3.5 -6){${Y_1(-1)}^{(1,2)}{Y_1(0)}^{(0,-1)}$}
\htext(3.5 -6.7){$\cdot {Y_2(-2)}^{(1,-1)}$}

\htext(9 -6){${Y_1(-1)}^{(1,6)}$}
\htext(9 -6.7){$\cdot {Y_2(-2)}^{(1,-2)}{Y_2(-1)}^{(0,-2)}$}
%%%%%%%%%%%%%%%%%%%%%%%%%%%%%%%%%%%%%%
\htext(-9 -9){${Y_1(-1)}^{(1,-3)}{Y_1(0)}^{(0,-3)}$}
\htext(-9 -9.7){$\cdot {Y_2(-2)}^{(1,0)}{Y_2(-1)}^{(0,3)}$}

\htext(-7.5 -11){${Y_1(-1)}^{(1,-2)}{Y_1(0)}^{(0,1)}$}
\htext(-7.5 -11.7){$\cdot {Y_2(-2)}^{(1,0)}{Y_2(-1)}^{(0,1)}
                       {Y_2(0)}^{(0,-1)}$}

\htext(-2.5 -9){${Y_1(-1)}^{(1,-1)}{Y_1(0)}^{(0,1)}$}
\htext(-2.5 -9.7){$\cdot {Y_1(1)}^{(0,-1)}{Y_2(-2)}^{(1,0)}$}

\htext(1 -11){${Y_1(-1)}^{(1,2)}{Y_1(0)}^{(0,2)}{Y_2(-2)}^{(1,-1)}$}
\htext(1 -11.7){$\cdot {Y_2(-1)}^{(0,-1)}{Y_2(0)}^{(0,-1)}$}

\htext(3.8 -9){${Y_1(-1)}^{(1,1)}{Y_1(0)}^{(0,-2)}$}
\htext(3.8 -9.7){$\cdot {Y_2(-2)}^{(1,-1)}{Y_2(-1)}^{(0,1)}$}

\htext(9 -11){${Y_1(-1)}^{(1,5)}{Y_1(0)}^{(0,-1)}$}
\htext(9 -11.7){$\cdot {Y_2(-2)}^{(1,-2)}{Y_2(-1)}^{(0,-1)}$}

\htext(10.5 -9){${Y_1(-1)}^{(1,9)}$}
\htext(10.5 -9.7){$\cdot {Y_2(-2)}^{(1,-3)}{Y_2(-1)}^{(0,-3)}$}
%%%%%%%%%%%%%%
\move(0 0)
\bsegment
\htext(-1.8 -1.5){$1$}
\htext(1.8 -1.5){$2$}
\esegment
\move(0 -2.7)
\bsegment
\htext(-5.7 -2.1){$1$}
\htext(-2.2 -2.1){$2$}
\htext(2.2 -2.1){$1$}
\htext(6.1 -2.1){$2$}
\esegment
\move(0 -6.2)
\bsegment
\htext(-9 -1.7){$1$}
\htext(-6 -1.7){$2$}
\htext(-2.9 -1.7){$1$}
\htext(0.2 -1.7){$2$}
\htext(3.1 -1.7){$1$}
\htext(5.6 -1.7){$2$}
\htext(7.5 -1.7){$1$}
\htext(9.8 -1.7){$2$}
\esegment
\move(0 -13.3)
%\drawbb
\end{texdraw}%
\caption{The monomial set $\M(\infty)$}\label{fig:1}
\end{figure}

Actually, as we will become apparent later, this set $\M(\infty)$ is
closed and connected under Kashiwara operators. Figure~\ref{fig:1} is
the top part of monomial set $\M(\infty)$.

We now introduce new expressions for elements of $\M(\infty)$.
First, we introduce the following notation.

\begin{definition}
For $u\in\Z_{\ge 0}$, $v\in\Z$, and $m\in\Z$, we use the notation
\begin{align}\label{order}
{X_j(m)}^{(u,v)}&=
\begin{cases}
{Y_j(m)}^{(u,v)}{Y_{j-1}(m+1)}^{(-u,-v)} &\text{for $j=1,2$},\\
{Y_1(m+1)}^{(2u,2v)}{Y_2(m+1)}^{(-u,-v)} &\text{for $j=3$},
\end{cases}\\
{X_0(m)}^{(u,v)}
&={Y_1(m+1)}^{(u,v)}{Y_1(m\!+\!2)}^{(-u,-v)},\\
{X_{\bar{j}}(m)}^{(u,v)}&=
\begin{cases}
{Y_{j-1}(m+(4-j))}^{(u,v)}{Y_j(m+(4-j))}^{(-u,-v)} &\text{for $j=1,2$},\\
{Y_2(m+1)}^{(u,v)}{Y_1(m+2)}^{(-2u,-2v)} &\text{for $j=3$}.
\end{cases}
\end{align}
Here, we set $Y_0(k)^{(u,v)}=1$.
\end{definition}

\begin{remark}
Using the above notation, we may write
\begin{align*}
A_1(m)&= X_1(m)^{(0,1)}X_2(m)^{(0,-1)}\\
      &= X_3(m-1)^{(0,1)}X_0(m-1)^{(0,-1)}\\
      &= X_0(m-1)^{(0,1)}X_{\bar3}(m-1)^{(0,-1)}\\
      &= X_{\bar2}(m-2)^{(0,1)}X_{\bar1}(m-2)^{(0,-1)},\\
A_2(m) &= X_2(m)^{(0,1)}X_3(m)^{(0,-1)}\\
       &= X_{\bar3}(m-1)^{(0,1)}X_{\bar2}(m-1)^{(0,-1)}.
\end{align*}
This is very useful when computing Kashiwara action on monomials
written in terms of $X_j(m)^{(u,v)}$ or $X_{\bar j}(m)^{(u,v)}$.
\end{remark}

\begin{prop}\label{propnk}
Consider elements of $\Me$ having the form
\begin{equation}~\label{equ:jin18}
\begin{aligned}
      M=& X_1(-1)^{(2,-b^{-1}_2-b^{-1}_3-b^{-1}_0-
                      b^{-1}_{\bar3}-b^{-1}_{\bar2}-b^{-1}_{\bar1})}
          X_2(-1)^{(0,b_2^{-1})}
          X_3(-1)^{(0,b_3^{-1})}\\
        & \cdot
          X_0(-1)^{(0,b_0^{-1})}
          X_{\bar3}(-1)^{(0,b_{\bar3}^{-1})}
          X_{\bar2}(-1)^{(0,b_{\bar2}^{-1})}
          X_{\bar1}(-1)^{(0,b_{\bar1}^{-1})}\\
        & \cdot
          X_2(-2)^{(1,-b_3^{-2})}
          X_3(-2)^{(0,b_3^{-2})}
\end{aligned}
\end{equation}
where $b_i^j\ge 0$ for all $i,j$ and $b_0^{-1}\le 1$.
Each element of $\M(\infty)$ may be written uniquely in
this form. Conversely, any element of this form is an element of
$\M(\infty)$.
\end{prop}
\begin{proof}
Given any monomial
\begin{equation*}
\begin{aligned}
      M=&{Y_1(-1)}^{(1,a_1^{-1})}{Y_1(0)}^{(0,a_1^0)}
        {Y_1(1)}^{(0,a_1^1)}{Y_1(2)}^{(0,a_1^2)}\\
        &\cdot {Y_2(-2)}^{(1,a_2^{-2})}{Y_2(-1)}^{(0,a_2^{-1})}
        {Y_2(0)}^{(0,a_2^0)}{Y_2(1)}^{(0,a_2^1)}\in\M(\infty),
\end{aligned}
\end{equation*}
through simple computation, we can obtain the expression
\begin{equation}~\label{equ:jin19}
\begin{aligned}
M=& X_1(-1)^{(2,a_1^{-1}+3a_2^{-2})}
    X_2(-1)^{(0,a_2^{-1}-a_2^{-2})}\\
  & \cdot X_3(-1)^{(0,t_3^{-1})}
    X_0(-1)^{(0,t_0^{-1})}
    X_{\bar 3}(-1)^{(0,t_{\bar3}^{-1})}\\
  & \cdot X_{\bar2}(-1)^{(0,-a_2^1)}
    X_{\bar1}(-1)^{(0,-a_1^2)}X_2(-2)^{(1,a_2^{-2})}
    X_3(-2)^{(0,-a_2^{-2})},
\end{aligned}
\end{equation}
where either
\begin{itemize}
\item $t_0^{-1}=0$,\ \
      $2t_3^{-1}=a_1^0+a_2^{-1}-a_2^{-2}$,\ \ and\
      $2t^{-1}_{\bar 3}=-a_1^1-a_2^1$,
\end{itemize}
or
\begin{itemize}
\item $t^{-1}_0=1$,\ \
      $2t^{-1}_3=a_1^0+a_2^{-1}-a_2^{-2}-1$,\ \ and\
      $2t^{-1}_{\bar 3}=-a_1^1-a_2^1-1$.
\end{itemize}
Since $M\in\M(\infty)$, from the conditions given in~\eqref{extend}, we
obtain the form given in~\eqref{equ:jin18}. Specifically, the element
$M_{\infty}={Y_1(-1)}^{(1,0)}{Y_2(-2)}^{(1,0)}$ corresponds to
${X_1(-1)}^{(2,0)}{X_2(-2)}^{(1,0)}$.

Conversely, given any monomial of the form \eqref{equ:jin18}, we
have
\begin{equation}\label{equ:jin7}
\begin{aligned}
M=&{Y_1(-1)}^{(1,-b^{-1}_2-b^{-1}_3-b^{-1}_0-b^{-1}_{\bar3}
                 -b^{-1}_{\bar2}-b^{-1}_{\bar1}+3b^{-2}_3)}\\
  &\cdot {Y_1(0)}^{(0,b^{-1}_0-b^{-1}_2+2b^{-1}_3)}
   {Y_1(1)}^{(0,-b^{-1}_0-2b^{-1}_{\bar 3}+b^{-1}_{\bar2})}
   {Y_1(2)}^{(0,-b^{-1}_{\bar1})}\\
  &\cdot {Y_2(-2)}^{(1,-b^{-2}_3)}
         {Y_2(-1)}^{(0,b^{-1}_2-b^{-2}_3)}
         {Y_2(0)}^{(0,b^{-1}_{\bar3}-b^{-1}_3)}
         {Y_2(1)}^{(0,-b^{-1}_{\bar2})}.
\end{aligned}
\end{equation}
It is now straightforward to check that $M\in\M(\infty)$. We have
thus shown that $\M(\infty)$ consists of elements of the
form~\eqref{equ:jin18}.

The uniqueness part may be proved through simple computation.
\end{proof}

\begin{remark}
There are other ways to write each element of $\M(\infty)$ as
products of the terms $X_j(m)^{(u,v)}$ and $X_{\bar
j}(m)^{(u,v)}$. The product form~\eqref{equ:jin18} was chosen
because it allows us to relate monomials of the set $\M(\infty)$
to tableaux in $\T(\infty)$ directly.
\end{remark}

Now, we translate the Kashiwara actions~\eqref{actions},
\eqref{actions2} into a form suitable for the new monomial expression
of $\M(\infty)$.

\begin{lemma}\label{anaysiskashi2}
The Kashiwara operator action on $\Me$ may be rewritten as given below
for elements
\begin{equation}~\label{equ:jin8}
\begin{aligned}
      M=& X_1(-1)^{(2,-b^{-1}_2-b^{-1}_3-b^{-1}_0-
                      b^{-1}_{\bar3}-b^{-1}_{\bar2}-b^{-1}_{\bar1})}
          X_2(-1)^{(0,b_2^{-1})}
          X_3(-1)^{(0,b_3^{-1})}\\
        & \cdot
          X_0(-1)^{(0,b_0^{-1})}
          X_{\bar3}(-1)^{(0,b_{\bar3}^{-1})}
          X_{\bar2}(-1)^{(0,b_{\bar2}^{-1})}
          X_{\bar1}(-1)^{(0,b_{\bar1}^{-1})}\\
        & \cdot
          X_2(-2)^{(1,-b_3^{-2})}
          X_3(-2)^{(0,b_3^{-2})}
\end{aligned}
\end{equation}
of $\M(\infty)$. Elements of the above form constitutes $\M(\infty)$
and this set is closed under Kashiwara operator actions.

\noindent \textup{$(1)$} Kashiwara actions $\tilde f_1$ and $\tilde
e_1$:
\begin{itemize}
\item Consider the following ordered sequence
      of some components of $M$.
      \begin{equation*}
      {X_{\bar1}(-1)}^{(0,b_{\bar1}^{-1})}
      {X_{\bar2}(-1)}^{(0,b_{\bar2}^{-1})}
      {X_{\bar3}(-1)}^{(0,b_{\bar3}^{-1})}
      {X_0(-1)}^{(0,b_0^{-1})}
      {X_3(-1)}^{(0,b_3^{-1})}
      {X_2(-1)}^{(0,b_2^{-1})}.
      \end{equation*}
\item Under each of the components
      \begin{equation*}
      {X_{\bar1}(-1)}^{(0,b_{\bar1}^{-1})},
      {X_0(-1)}^{(0,b_0^{-1})},
      {X_2(-1)}^{(0,b_2^{-1})},
      \end{equation*}
      given in the above sequence, write $b_j^{-1}$-many $1$'s and
      under
      ${X_{\bar3}(-1)}^{(0,b_{\bar3}^{-1})}$,
      write $(2b_{\bar 3}^{-1})$-many $1$'s.
      Also, under each of the components
      \begin{equation*}
      {X_{\bar2}(-1)}^{(0,b_{\bar2}^{-1})},{X_0(-1)}^{(0,b_0^{-1})},
      \end{equation*}
      write $b_j^{-1}$-many $0$'s and under
      ${X_3(-1)}^{(0,b_3^{-1})}$, write $(2b_3^{-1})$-many $0$'s.
\item From this sequence of $1$'s and $0$'s, successively
      cancel out each $(0,1)$-pair to obtain a sequence of $1$'s
      followed by $0$'s \textup{(}reading from left to right\textup{)}.
      This remaining $1$ and $0$ sequence is called
      the \defi{$1$-signature of $M$}.
\item Depending on the component $X$ corresponding
      to the leftmost $0$ of the $1$-signature of $M$,
      we define $\tilde f_1 M$ as follows\,:
      \begin{equation}\label{lem1}
      \tilde f_1 M=
      \begin{cases}
      MX_{\bar2}(-1)^{(0,-1)}X_{\bar1}(-1)^{(0,1)}=MA_1(1)^{-1}
            \quad\textup{if $X={X_{\bar2}(-1)}^{(0,b_{\bar2}^{-1})}$},\\
      MX_0(-1)^{(0,-1)}X_{\bar3}(-1)^{(0,1)}=MA_1(0)^{-1}
            \quad\textup{if $X={X_0(-1)}^{(0,b_0^{-1})}$},\\
      MX_3(-1)^{(0,-1)}X_0(-1)^{(0,1)}=MA_1(0)^{-1}
            \quad\textup{if $X={X_3(-1)}^{(0,b_3^{-1})}$}.
      \end{cases}
      \end{equation}
      We define
      \begin{equation}\label{lem3}
      \tilde f_1 M=M{X_1(-1)}^{(0,-1)}{X_{2}(-1)}^{(0,1)}
      =MA_1(-1)^{-1}
      \end{equation}
      if no $0$ remains.
\item Depending on the component $X$ corresponding
      to the rightmost $1$ of the $1$-signature of $M$,
      we define $\tilde e_1 M$ as follows\,:
      \begin{equation*}
      \tilde e_1 M=
      \begin{cases}
      M X_{\bar2}(-1)^{(0,1)}X_{\bar1}(-1)^{(0,-1)}=M A_1(1)
            \quad \textup{if $X={X_{\bar 1}(-1)}^{(0,b_{\bar 1}^{-1})}$},\\
      M X_0(-1)^{(0,1)}X_{\bar3}(-1)^{(0,-1)}=M A_1(0)
            \quad \textup{if $X={X_{\bar3}(-1)}^{(0,b_{\bar3}^{-1})}$},\\
      M X_3(-1)^{(0,1)}X_0(-1)^{(0,-1)}=M A_1(0)
            \quad \textup{if $X={X_0(-1)}^{(0,b_0^{-1})}$},\\
      M {X_1(-1)}^{(0,1)}{X_{2}(-1)}^{(0,-1)}=M A_1(-1)
            \ \textup{if $X={X_2(-1)}^{(0,b_2^{-1})}$}.
      \end{cases}
      \end{equation*}
      We define $\tilde e_1 M=0$ if no $1$ remains.
\end{itemize}
\textup{$(2)$} Kashiwara actions $\tilde f_2$ and $\tilde e_2$ :
\begin{itemize}
\item Consider the following finite ordered sequence
      of some components of $M$.
      \begin{equation*}
      {X_{\bar2}(-1)}^{(0,b_{\bar2}^{-1})}
      {X_{\bar3}(-1)}^{(0,b_{\bar3}^{-1})}
      {X_3(-1)}^{(0,b_3^{-1})}
      {X_2(-1)}^{(0,b_2^{-1})}
      {X_3(-2)}^{(0,b_3^{-2})}.
      \end{equation*}
\item Under each of the components
      \begin{equation*}
      {X_{\bar2}(-1)}^{(0,b_{\bar2}^{-1})},
      {X_3(-1)}^{(0,b_3^{-1})},
      {X_3(-2)}^{(0,b_3^{-2})},
      \end{equation*}
      from the above sequence, write $b_j^k$-many $1$'s,
      and under each
      \begin{equation*}
      {X_{\bar3}(-1)}^{(0,b_{\bar3}^{-1})},
      {X_2(-1)}^{(0,b_2^{-1})},
      \end{equation*}
      write $b_j^{-1}$-many $0$'s.
\item From this sequence of $1$'s and $0$'s, successively
      cancel out each $(0,1)$-pair to obtain a sequence of $1$'s
      followed by $0$'s.
      This remaining $1$ and $0$ sequence is called
      the \defi{$2$-signature of $M$}.
\item Depending on the component $X$
      corresponding to the leftmost $0$ of the $2$-signature of $M$,
      we define $\tilde f_2 M$ as follows\,:
      \begin{equation}\label{lem2}
      \tilde f_2 M=
      \begin{cases}
      M X_{\bar3}(-1)^{(0,-1)}X_{\bar2}(-1)^{(0,1)}=M A_2(0)^{-1}
            \ \ \ \textup{if $X={X_{\bar3}(-1)}^{(0,b_{\bar3}^{-1})}$},\\
      MX_2(-1)^{(0,-1)}X_{3}(-1)^{(0,1)}=MA_2(-1)^{-1}
             \ \textup{if $X={X_2(-1)}^{(0,b_2^{-1})}$}.
      \end{cases}
      \end{equation}
      We define
      \begin{equation}\label{lem4}
      \tilde f_2 M=M{X_2(-2)}^{(0,-1)}{X_{3}(-2)}^{(0,1)}
      =MA_2(-2)^{-1}
      \end{equation}
      if no $0$ remains.
\item Depending on the component $X$ corresponding
      to the rightmost $1$ of the $2$-signature of $M$,
      we define $\tilde e_2 M$ as follows\,:
      \begin{equation*}
      \tilde e_2 M=
      \begin{cases}
      M X_{\bar3}(-1)^{(0,1)}X_{\bar2}(-1)^{(0,-1)}=MA_2(0)
            \quad\ \ \textup{if $X={X_{\bar2}(-1)}^{(0,b_{\bar2}^{-1})}$},\\
      M X_2(-1)^{(0,1)}X_{3}(-1)^{(0,-1)}=MA_2(-1)
            \quad \textup{if $X={X_{3}(-1)}^{(0,b_{3}^{-1})}$},\\
      M {X_2(-2)}^{(0,1)}{X_{3}(-2)}^{(0,-1)}=MA_2(-2)
            \quad \textup{if $X={X_3(-2)}^{(0,b_3^{-2})}$}.
      \end{cases}
      \end{equation*}
      We define $\tilde e_2 M=0$ if no $1$ remains.
\end{itemize}
\end{lemma}
\begin{proof}
We first show that the action of these operators is closed on
$\M(\infty)$.

For $M\in\M(\infty)$, if the $i$-signature of $M$ contains at least one
$0$, then the exponent of component ${X_j(-1)}^{(0,b_j^{-1})}$
corresponding to the left-most $0$ shows the property $(0,b_j^{-1})\ge
(0,1)$. In particular, in dealing with the action of $\tilde f_1$, if
the left-most $0$ corresponds to ${X_3(-1)}^{(0,b_3^{-1})}$, it means
that the exponent $b_0^{-1}$ of ${X_0(-1)}^{(0,b_0^{-1})}$, a component
of $M$, is $0$, due to the $(0,1)$-pair cancellation rule. Thus the
monomial $\fit M$ defined in~\eqref{lem1} and \eqref{lem2} is contained
in $\M(\infty)$. In the case where $i$-signature of $M$ contains no
$0$, the exponents of the components $X_i(-i)$ of $M$ show the property
$\ge (0,1)$. Thus $\fit M$ given in~\eqref{lem3} and \eqref{lem4} also
are in $\M(\infty)$. So the set $\M(\infty)$ is closed under the above
operator $\fit$.

As we can see in equations~\eqref{lem1} to~\eqref{lem4}, for each
$M\in\M(\infty)$, $\fit M$ can also be expressed in form
$MA_i(m)^{-1}$. To show that this operation is just another
interpretation of the Kashiwara operator $\fit$ given on $\Me$,
restricted to $\M(\infty)$, it is enough to show that $m_f$ defined
in~\eqref{emf} for each $M$ is equal to $m$ of $MA_i(m)^{-1}$ given in
equations~\eqref{lem1} to~\eqref{lem4}.

Given a monomial $M\in\M(\infty)$, we can express it in the
following two forms.
\begin{align}
      M=& X_1(-1)^{(2,-b^{-1}_2-b^{-1}_3-b^{-1}_0-
                      b^{-1}_{\bar3}-b^{-1}_{\bar2}-b^{-1}_{\bar1})}
          X_2(-1)^{(0,b_2^{-1})}
          X_3(-1)^{(0,b_3^{-1})}\label{eq22}\\
        & \cdot
          X_0(-1)^{(0,b_0^{-1})}
          X_{\bar3}(-1)^{(0,b_{\bar3}^{-1})}
          X_{\bar2}(-1)^{(0,b_{\bar2}^{-1})}
          X_{\bar1}(-1)^{(0,b_{\bar1}^{-1})}\notag\\
        & \cdot
          X_2(-2)^{(1,-b_3^{-2})}
          X_3(-2)^{(0,b_3^{-2})}\notag\\
       =&{Y_1(-1)}^{(1,-b^{-1}_2-b^{-1}_3-b^{-1}_0-b^{-1}_{\bar3}
                 -b^{-1}_{\bar2}-b^{-1}_{\bar1}+3b^{-2}_3)}\label{eq23}\\
  &\cdot {Y_1(0)}^{(0,b^{-1}_0-b^{-1}_2+2b^{-1}_3)}
   {Y_1(1)}^{(0,-b^{-1}_0-2b^{-1}_{\bar 3}+b^{-1}_{\bar2})}
   {Y_1(2)}^{(0,-b^{-1}_{\bar1})}\notag\\
  &\cdot {Y_2(-2)}^{(1,-b^{-2}_3)}
         {Y_2(-1)}^{(0,b^{-1}_2-b^{-2}_3)}
         {Y_2(0)}^{(0,b^{-1}_{\bar3}-b^{-1}_3)}
         {Y_2(1)}^{(0,-b^{-1}_{\bar2})}\notag.
\end{align}
If the $1$-signature of $M$ contains at least one $0$ and $X$ is
the component corresponding to the left-most $0$ in the
$1$-signature of $M$, then we can obtain
\begin{equation}
m_f=\textup{min}\{j\in\Z\ \vert\
    \textup{max}\{ \sum_{k\le j} y_1(k)\}\}=
\begin{cases}
1 &\text{if $X={X_{\bar 2}(-1)}^{(0,b_{\bar 2}^{-1})}$},\\
0 &\text{if $X={X_0(-1)}^{(0,b_0^{-1})}$},\\
0 &\text{if $X={X_3(-1)}^{(0,b_3^{-1})}$},
\end{cases}
\end{equation}
where $y_1(k)$ is the exponent of $Y_1(k)$ appearing in $M$ given by
expression~\eqref{eq23}. If the $2$-signature of $M$ contains at least
one $0$ and $X$ is the component corresponding to the left-most $0$ in
the $2$-signature of $M$, then we can obtain
\begin{equation}
m_f=\textup{min}\{j\in\Z\ \vert\
    \textup{max}\{ \sum_{k\le j} y_2(k)\}\}=
\begin{cases}
0 &\text{for $X={X_{\bar 3}(-1)}^{(0,b_{\bar 3}^{-1})}$},\\
-1 &\text{for $X={X_2(-1)}^{(0,b_2^{-1})}$},
\end{cases}
\end{equation}
where $y_2(k)$ is the exponent of $Y_2(k)$ appearing in $M$ given by
expression~\eqref{eq23}.

If the $i$-signature of $M$ contains no $0$,
\begin{equation}
m_f=\textup{min}\{j\in\Z\ \vert\
     \textup{max}\{ \sum_{k\le j} y_i(k)\}\}=-i,
\end{equation}
where $y_i(k)$ is the exponent of $Y_i(k)$ appearing in $M$ given by
expression~\eqref{eq23}.

In all cases, we can confirm that $m_f=m$, where $m$ is given through
equations~\eqref{lem1} to~\eqref{lem4} stating $\fit M = MA_i(m)^{-1}$.
Proof for the statements concerning $\eit$ may be done in a similar
manner.
\end{proof}

In the above lemma, we showed that the action of Kashiwara
operators~\eqref{actions} and~\eqref{actions2} on $\Me$ satisfy the
following properties\;:
\begin{equation}
\fit \M(\infty) \subset \M(\infty), \qquad \eit \M(\infty) \subset
\M(\infty) \cup \{0\} \qquad \text{for all} \ \ i\in I.
\end{equation}
Thus we obtain the following result.

\begin{prop}\label{zaq}
The set $\M(\infty)$ forms a $\uq(G_2)$-subcrystal of $\Me$.
\end{prop}

\begin{figure}
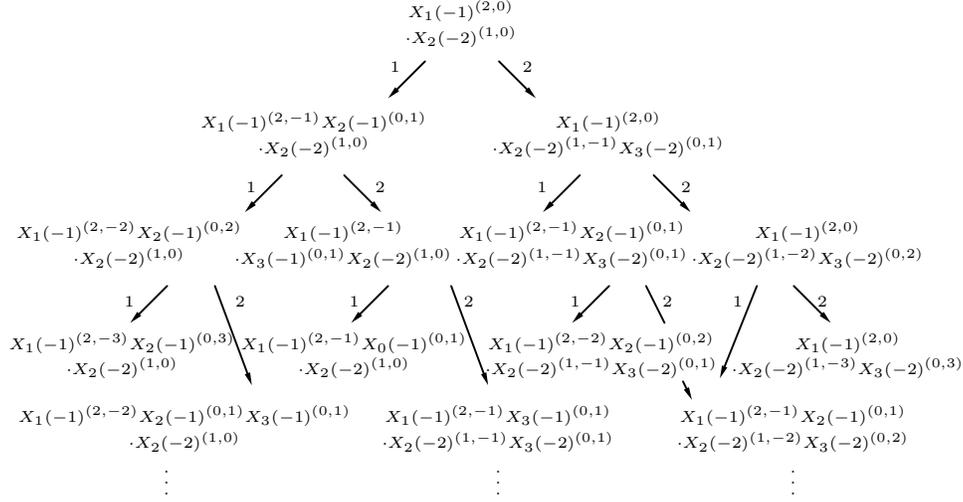

\centering
\begin{texdraw}%
\drawdim in \arrowheadsize l:0.065 w:0.03 \arrowheadtype t:F
\fontsize{6}{6}\selectfont \textref h:C v:C \drawdim em
\setunitscale 1.9
\move(-1 -1.4)\ravec(-1 -1)%%%%%%%%%%%%%%%%%
\move(1 -1.4)\ravec(1 -1)%%%%%%%%%%%%%%%%%
\move(-3.2 -4.5)\ravec(1 -1)%%%%%%%%%%%%%%%%%
\move(-4.9 -4.5)\ravec(-1 -1)%%%%%%%%%%%%%%%%%
\move(3.2 -4.5)\ravec(-1 -1)%%%%%%%%%%%%%%%%%
\move(5.2 -4.5)\ravec(1 -1)%%%%%%%%%%%%%%%%%
\move(-8 -7.5)\ravec(-1 -1)%%%%%%%%%%%%%%%%%
\move(-6.7 -7.5)\ravec(1 -2.7)%%%%%%%%%%%%%%%%%
\move(-2 -7.5)\ravec(-1 -1)%%%%%%%%%%%%%%%%%
\move(-0.3 -7.5)\ravec(1 -2.7)%%%%%%%%%%%%%%%%%
\move(4 -7.5)\ravec(-1 -1)%%%%%%%%%%%%%%%%%
\move(5 -7.5)\rlvec(0.5 -1)%%%%%%%%%%%%%%%%%
\move(6 -10.1)\ravec(0.25 -0.5)%%%%%%%%%%%%%%%%%
\move(8 -7.5)\ravec(-1 -2.5)%%%%%%%%%%%%%%%%%
\move(9 -7.5)\ravec(1 -1)%%%%%%%%%%%%%%%%%
\htext(-8 -12.6){$\vdots$}%%%%%%%%%%%%%%%%%
\htext(9 -12.6){$\vdots$}%%%%%%%%%%%%%%%%%
\htext(1 -12.6){$\vdots$}%%%%%%%%%%%%%%%%%
%%%%%%%%%%%%%%%%%%%%%%%%%%%%%%%%%%%%%%%%%%%%%%%
\htext(0 0){${X_1(-1)}^{(2,0)}$} \htext(0 -0.7){$\cdot
{X_2(-2)}^{(1,0)}$}
%%%%%%%%%%%%%%%%%%%%%%%%%%%%%%%%%%%%%%%%%%%%%%%
\htext(-4 -3){${X_1(-1)}^{(2,-1)}{X_2(-1)}^{(0,1)}$}
\htext(-4 -3.7){$\cdot {X_2(-2)}^{(1,0)}$}

\htext(4 -3){${X_1(-1)}^{(2,0)}$}
\htext(4 -3.7){$\cdot {X_2(-2)}^{(1,-1)}{X_3(-2)}^{(0,1)}$}
%%%%%%%%%%%%%%%%%%%%%%%%%%%%%%%%%%%%%%%%%%%%%%
\htext(-9 -6){${X_1(-1)}^{(2,-2)}{X_2(-1)}^{(0,2)}$}
\htext(-9 -6.7){$\cdot {X_2(-2)}^{(1,0)}$}

\htext(-3.2 -6){${X_1(-1)}^{(2,-1)}$}
\htext(-3.2 -6.7){$\cdot {X_3(-1)}^{(0,1)}{X_2(-2)}^{(1,0)}$}

\htext(3 -6){${X_1(-1)}^{(2,-1)}{X_2(-1)}^{(0,1)}$}
\htext(3 -6.7){$\cdot {X_2(-2)}^{(1,-1)}{X_3(-2)}^{(0,1)}$}

\htext(9.4 -6){${X_1(-1)}^{(2,0)}$}
\htext(9.4 -6.7){$\cdot {X_2(-2)}^{(1,-2)}{X_3(-2)}^{(0,2)}$}
%%%%%%%%%%%%%%%%%%%%%%%%%%%%%%%%%%%%%%%%%%%%%%%%%%%%%%%%%%%%%%%%%%%%
\htext(-9.2 -9){${X_1(-1)}^{(2,-3)}{X_2(-1)}^{(0,3)}$}
\htext(-9.2 -9.7){$\cdot {X_2(-2)}^{(1,0)}$}

\htext(-7.5 -11){${X_1(-1)}^{(2,-2)}{X_2(-1)}^{(0,1)}{X_3(-1)}^{(0,1)}$}
\htext(-7.5 -11.7){$\cdot {X_2(-2)}^{(1,0)}$}

\htext(-2.9 -9){${X_1(-1)}^{(2,-1)}{X_0(-1)}^{(0,1)}$}
\htext(-2.9 -9.7){$\cdot {X_2(-2)}^{(1,0)}$}

\htext(1 -11){${X_1(-1)}^{(2,-1)}{X_3(-1)}^{(0,1)}$}
\htext(1 -11.7){$\cdot {X_2(-2)}^{(1,-1)}{X_3(-2)}^{(0,1)}$}

\htext(3.8 -9){${X_1(-1)}^{(2,-2)}{X_2(-1)}^{(0,2)}$}
\htext(3.8 -9.7){$\cdot {X_2(-2)}^{(1,-1)}{X_3(-2)}^{(0,1)}$}

\htext(9 -11){${X_1(-1)}^{(2,-1)}{X_2(-1)}^{(0,1)}$}
\htext(9 -11.7){$\cdot {X_2(-2)}^{(1,-2)}{X_3(-2)}^{(0,2)}$}

\htext(10.5 -9){${X_1(-1)}^{(2,0)}$}
\htext(10.5 -9.7){$\cdot {X_2(-2)}^{(1,-3)}{X_3(-2)}^{(0,3)}$}
%%%%%%%%%%%%%%
\move(0 0) \bsegment \htext(-1.8 -1.5){$1$} \htext(1.8 -1.5){$2$}
\esegment \move(0 -2.7) \bsegment \htext(-5.7 -2.1){$1$}
\htext(-2.2 -2.1){$2$} \htext(2.2 -2.1){$1$} \htext(6.1 -2.1){$2$}
\esegment \move(0 -6.2) \bsegment \htext(-9 -1.7){$1$} \htext(-6
-1.7){$2$} \htext(-2.9 -1.7){$1$} \htext(0.2 -1.7){$2$} \htext(3.1
-1.7){$1$} \htext(5.6 -1.7){$2$} \htext(7.5 -1.7){$1$} \htext(9.8
-1.7){$2$} \esegment \move(0 -13.3)
\end{texdraw}
\caption{Crystal $\M(\infty)$ for type $G_2$}\label{fig:3}
\end{figure}

Figures~\ref{fig:3} illustrates the top part of crystal
$\M(\infty)$ for finite type $G_2$. It was obtained by applying
the Kashiwara actions introduced in Lemma~\ref{anaysiskashi2} on
the new expression for elements of $\M(\infty)$. Readers may want
to compare this with Figure~\ref{fig:1}.

Actually, from property of the crystal structure of $\Me$ we can
obtain more general results.

\begin{definition}
Fix any set of positive integers $p_i$ and any integer $r$. Consider
elements of $\Me$ having the form
\begin{equation}\label{extend2}
\begin{aligned}
      M=&{Y_1(r-1)}^{(p_1,a_1^{-1})}{Y_1(r)}^{(0,a_1^0)}
        {Y_1(r+1)}^{(0,a_1^1)}{Y_1(r+2)}^{(0,a_1^2)}\\
        &\cdot {Y_2(r-2)}^{(p_2,a_2^{-2})}{Y_2(r-1)}^{(0,a_2^{-1})}
        {Y_2(r)}^{(0,a_2^0)}{Y_2(r+1)}^{(0,a_2^1)}
\end{aligned}
\end{equation}
satisfying the same condition given to~\eqref{extend}. When $a_i^j=0$
for all $i,j$, this reduces to
\begin{equation}\label{hi1}
      M={Y_1(r-1)}^{(p_1,0)}{Y_2(r-2)}^{(p_2,0)}.
\end{equation}
      We denote by $\M(p_1,p_2;r;\infty)$
      the set of all monomials of this form
      and write the monomial of~\eqref{hi1} as $M_{(p_1,p_2;r;\infty)}$.
\end{definition}
A result similar to Proposition~\ref{propnk} may be obtained for
$M_{(p_1,p_2;r;\infty)}$.
\begin{prop}
Each element of $\M(p_1,p_2;r;\infty)$ may be written uniquely in the
form
\begin{equation}~\label{equ:jin6}
\begin{aligned}
      M=& X_1(r-1)^{(p_1+p_2,-b^1_2-b^1_3-b^1_0-
                      b^1_{\bar3}-b^1_{\bar2}-b^1_{\bar1})}
          X_2(r-1)^{(0,b_2^1)}
          X_3(r-1)^{(0,b_3^1)}\\
        & \cdot
          X_0(r-1)^{(0,b_0^1)}
          X_{\bar3}(r-1)^{(0,b_{\bar3}^1)}
          X_{\bar2}(r-1)^{(0,b_{\bar2}^1)}
          X_{\bar1}(r-1)^{(0,b_{\bar1}^1)}\\
        & \cdot
          X_2(r-2)^{(p_2,-b_3^2)}
          X_3(r-2)^{(0,b_3^2)}
\end{aligned}
\end{equation}
where $b_i^j\ge 0$ for all $i,j$ and $b_0^1\le 1$.
Conversely, any element in $\Me$ of this form is
an element of $\M(p_1,p_2;r;\infty)$.
\end{prop}

We believe the readers can easily write down the process for
\emph{change of variable} similar to that given
by~\eqref{equ:jin19} and~\eqref{equ:jin7} for
$\M(p_1,p_2;r;\infty)$.

The set $\M(\infty)$ is a special case of this set
$\M(p_1,p_2;r;\infty)$ corresponding to $r=0$ and $p_i=1$ for all $i\in
I$.

\begin{remark}\label{lem6}
It is possible to obtain the result of Lemma~\ref{anaysiskashi2} also
for the case $\M(p_1,p_2;r;\infty)$. Thus we can state that the set
$\M(p_1,p_2;r;\infty)$ forms a $\uq(G_2)$-subcrystal of $\Me$.
\end{remark}

\begin{prop}\label{cor4}
The set $\M(p_1,p_2;r;\infty)$ forms a subcrystal of $\Me$ isomorphic
to $\M(\infty)$ as a $\uq(G_2)$-crystal.
\end{prop}
\begin{proof}
As mentioned in Remark~\ref{lem6}, we can show that the set
$\M(p_1,p_2;r;\infty)$ forms a $\uq(G_2)$-subcrystal of $\Me$. Let us
show that the crystal $\M(p_1,p_2;r;\infty)$ is isomorphic to
$\M(\infty)$ as a $\uq(G_2)$-crystal.

First, we define a canonical map
$\phi:\M(\infty)\rightarrow\M(p_1,p_2;r;\infty)$ by setting
\begin{equation}
\begin{aligned}
\phi(M)=& X_1(r-1)^{(p_1+p_2,-b^1_2-b^1_3-b^1_0-
                      b^1_{\bar3}-b^1_{\bar2}-b^1_{\bar1})}
          X_2(r-1)^{(0,b_2^1)}
          X_3(r-1)^{(0,b_3^1)}\\
        & \cdot
          X_0(r-1)^{(0,b_0^1)}
          X_{\bar3}(r-1)^{(0,b_{\bar3}^1)}
          X_{\bar2}(r-1)^{(0,b_{\bar2}^1)}
          X_{\bar1}(r-1)^{(0,b_{\bar1}^1)}\\
        & \cdot
          X_2(r-2)^{(p_2,-b_3^2)}
          X_3(r-2)^{(0,b_3^2)}
\end{aligned}
\end{equation}
for $M$ of the form~\eqref{equ:jin18}. Note that the monomial
$M_{\infty}$ of $\M(\infty)$ is mapped onto the vector
$M_{(p_1,p_2;r;\infty)}$. It is obvious that this map $\phi$ is
well-defined and that it is actually bijective.

The outputs of the functions $\wt$, $\vphi_i$, and $\veps_i$, defined
in~\eqref{structure1}, \eqref{structure2}, and \eqref{structure3} do
not depend on $r$ or on any fixed positive integers $p_i$. Using
Lemma~\ref{anaysiskashi2} and its counterpart for
$\M(p_1,p_2;r;\infty)$, we can easily show that the map $\phi$ commutes
with the Kashiwara operators. Hence, the set $\M(p_1,p_2;r;\infty)$
forms a subcrystal of $\Me$ isomorphic to $\M(\infty)$ as a
$\uq(G_2)$-crystal.
\end{proof}

\begin{remark}
The monomial of~\eqref{extend2} is the element of
$\M(p_1,p_2;r;\infty)$ corresponding to the monomial~\eqref{extend} of
$\M(\infty)$ under the natural isomorphism $\phi$ mentioned in the
proof of Proposition~\ref{cor4}.
\end{remark}

\begin{remark}\label{wkwmd3}
Since $\vepst_i(M_{\infty})=(0,0)$ for all $i\in I$, we have
$\eit(M_{\infty})=0$ for all $i\in I$. And
$\wtt(M_{\infty})=\sum_{i\in I} (1,0)\La_i$, so we have
$\wt(M_{\infty})$ $=0$. Note that $M_{\infty}$ satisfies the
condition for a monomial $M$ mentioned in the conjecture
introduced in Section~\ref{Nakajima monomials}.

The monomial $M_{(p_1,p_2;r;\infty)}$ with
$\wtt(M_{(p_1,p_2;r;\infty)})=\sum_{i\in I} (p_i,0)\La_i$, also
satisfies the condition for a monomial $M$ mentioned in the conjecture.
\end{remark}

\begin{remark}
It should be clear from the proof of Proposition~\ref{cor4}, that in
developing any theory for ${\M(p_1,p_2;r;\infty)}$ the actual
values of integer $r$ or $(p_1,p_2)$ will not be very important.
Arguments made for any set of such values can easily be adapted to
applied to other set of such values. Hence, we shall concentrate on the
theory for $\M(\infty)$ only.
\end{remark}

Now, we will show that $\M(\infty)$ is a new realization of
$\hwc(\infty)$ by giving a crystal isomorphism. Recall from
Theorem~\ref{HL} that the set $\T(\infty)$ gives a realization of the
crystal $\hwc(\infty)$.

Here is one of our two main realization theorems.

\begin{theorem}\label{mainreal5}
There exists a $\uq(G_2)$-crystal isomorphism
\begin{equation}
\B(\infty) \overset{\sim} \longrightarrow \T(\infty)
\overset {\sim} \longrightarrow \M(\infty)
\end{equation}
which maps $T_\infty$ to $M_{\infty}$.
\end{theorem}
\begin{proof}
We define a canonical map $\Theta:\T(\infty)\rightarrow\M(\infty)$
by setting, for each tableau $T\in\T(\infty)$ with second row
consists of $b^2_3$-many $3$-boxes and just one $2$-box, and with
first row consists of $b^1_j$-many $j$-boxes, for each $j\succ 1$,
and $(b^2_3+2)$-many $1$-boxes, $\Theta(T)=M$, where
\begin{equation*}
\begin{aligned}
      M=& X_1(-1)^{(2,-b^1_2-b^1_3-b^1_0-
                      b^1_{\bar3}-b^1_{\bar2}-b^1_{\bar1})}
          X_2(-1)^{(0,b_2^1)}
          X_3(-1)^{(0,b_3^1)}\\
        & \cdot
          X_0(-1)^{(0,b_0^1)}
          X_{\bar3}(-1)^{(0,b_{\bar3}^1)}
          X_{\bar2}(-1)^{(0,b_{\bar2}^1)}
          X_{\bar1}(-1)^{(0,b_{\bar1}^1)}\\
        & \cdot
          X_2(-2)^{(1,-b_3^2)}
          X_3(-2)^{(0,b_3^2)}\in\M(\infty).
\end{aligned}
\end{equation*}
It is obvious that this map $\Theta$ is well-defined and that it is
actually bijective.

The action of Kashiwara operators on $\M(\infty)$ given in
Lemma~\ref{anaysiskashi2} follows the process for defining it on
$\T(\infty)$. Hence, the map $\Theta$ naturally commutes with the
Kashiwara operators $\fit$ and $\eit$. Other parts of the proof are
similar or easy.
\end{proof}

\begin{remark}\label{final2}
 From Proposition~\ref{cor4} and Theorem~\ref{mainreal5}, we can conclude
that the crystal $\M(p_1,p_2;r;\infty)$ is also
$\uq(G_2)$-crystal isomorphic to $\B(\infty)$.
\end{remark}

\begin{example}
We illustrate the correspondence between $\T(\infty)$ and $\M(\infty)$.
A monomial of $\M(\infty)$
\begin{align*}
M=& {Y_1(-1)}^{(1,1)}{Y_1(1)}^{(0,-5)}{Y_1(2)}^{(0,-1)}\\
  &\cdot {Y_2(-2)}^{(1,-2)}{Y_2(-1)}^{(0,-1)}{Y_2(0)}^{(0,2)}
\end{align*}
can be expressed as
\begin{align*}
M=& {X_1(-1)}^{(2,-5)}
    {X_2(-1)}^{(0,1)}
    {X_0(-1)}^{(0,1)}
    {X_{\bar 3}(-1)}^{(0,2)}
    {X_{\bar 1}(-1)}^{(0,1)}\\
  &\cdot {X_2(-2)}^{(1,-2)}
         {X_3(-2)}^{(0,2)}
\end{align*}
by~(\ref{equ:jin19}). Hence we have the following marginally large
tableau as the image of $M$ under $\Theta^{-1}$.
\begin{equation*}
\Theta^{-1}(M)=
\raisebox{-1.3em}{
\begin{texdraw}%
\drawdim in
\arrowheadsize l:0.065 w:0.03
\arrowheadtype t:F
\fontsize{5}{5}\selectfont
\textref h:C v:C
\drawdim em
\setunitscale 1.6
\move(0 2)\rlvec(9 0)
\move(0 1)\rlvec(9 0)\rlvec(0 1)
\move(0 0)\rlvec(3 0)\rlvec(0 2)
\move(0 2)\rlvec(0 -2)
\move(1 0)\rlvec(0 2)
\move(2 0)\rlvec(0 2)
\move(4 1)\rlvec(0 1)
\move(5 1)\rlvec(0 1)
\move(6 1)\rlvec(0 1)
\move(7 1)\rlvec(0 1)
\move(8 1)\rlvec(0 1)
\htext(0.5 1.5){$1$}
\htext(1.5 1.5){$1$}
\htext(2.5 1.5){$1$}
\htext(3.5 1.5){$1$}
\htext(4.5 1.5){$2$}
\htext(5.5 1.5){$0$}
\htext(6.5 1.5){$\bar 3$}
\htext(7.5 1.5){$\bar 3$}
\htext(8.5 1.5){$\bar 1$}
\htext(0.5 0.5){$2$}
\htext(1.5 0.5){$3$}
\htext(2.5 0.5){$3$}
\end{texdraw}}%
\in\T(\infty).
\end{equation*}
\end{example}

\begin{remark}
Note that Remark~\ref{wkwmd3} and Remark~\ref{final2} provide
evidence supporting the conjecture introduced in
Section~\ref{Nakajima monomials}, for the $G_2$ case.
\end{remark}

%%%%%%%%%%%%%%%%%%%%%%%%%%%%%%%%%%%%%%%%%%%%%%%%%%%%%%%%%%%%%%%%%%%%%%%%%%
\section{Cliff's description of $\B(\infty)$}%
\label{Cliff's description}

Let us recall the abstract crystal $\B_i=\{b_i(k)\vert k\in{\bf Z}\}$
introduced in~\cite{MR1240605} for each $i\in I$. It has the following
maps defining the crystal structure.
\begin{alignat*}{3}
\wt b_i(k)      &=k\alpha_i,\\
\vphi_i(b_i(k)) &=k,         & \veps_i(b_i(k))&=-k,\\
\vphi_i(b_j(k)) &=-\infty,   & \veps_i(b_j(k))&=-\infty,
& &\text{ for }i\ne j,\\
\fit(b_i(k))    &=b_i(k-1),  \quad & \eit(b_i(k))&=b_i(k+1),\\
\fit(b_j(k))    &=0,         & \eit(b_j(k))&=0, & &\text{ for
}i\ne j.
\end{alignat*}
 From now on, we will denote the element $b_i(0)$ by $b_i$.
We next cite the tensor product rule on crystals.

\begin{prop}\textup{(\cite{MR1240605})}~\label{propq}
Let $\B^k(1\le k\le n)$ be crystals with $b^k\in \B^k$. We set
\begin{equation}\label{kashi}
a_k=\veps_i(b^k)-\sum_{1\le v<k} \langle h_i,\wt(b^v)\rangle.
\end{equation} Then we have
\begin{enumerate}
\item $\eit(b^1\otimes\cdots\otimes
      b^n)=b^1\otimes\cdots\otimes b^{k-1}\otimes
      \eit b^k\otimes b^{k+1}\otimes\cdots\otimes b^n$\\
      \mbox{}\hspace{20mm} if $a_k>a_v$ for $1\le v<k$
      and $a_k\ge a_v$ for $k<v\le n$,
\item $\fit(b^1\otimes\cdots\otimes
      b^n)=b^1\otimes\cdots\otimes b^{k-1}\otimes
      \fit b^k\otimes b^{k+1}\otimes\cdots\otimes b^n$\\
      \mbox{}\hspace{20mm} if $a_k\ge a_v$ for $1\le v<k$
      and $a_k>a_v$ for $k<v\le n$.
\end{enumerate}
\end{prop}

Kashiwara has shown~\cite{MR1240605} the existence of an injective
strict crystal morphism
\begin{equation}
\Psi:\B(\infty)\to\B(\infty)\otimes
     \B_{i_k}\otimes\B_{i_{k-1}}\otimes\cdots\otimes\B_{i_1}
\end{equation}
which sends the highest weight element $u_\infty$ to $u_\infty\otimes
b_{i_k}\otimes\cdots\otimes b_{i_1}$, for any sequence
$S={i_1,i_2,\cdots,i_k}$ of numbers in the index set $I$ of simple
roots. In~\cite{MR1614241}, Cliff uses this to give a combinatorial
description of $\B(\infty)$ for all finite classical types, with a
specific choice of sequence $S$. It is our goal to do this for type
$G_2$.

\begin{prop}
We define
\begin{align*}
\B(1)= \B_1\otimes
        \B_2\otimes
        \B_1\otimes
        \B_2\otimes
        \B_1 \ \text{and}\ \
\B(2)= \B_2.
\end{align*}
Consider the subset of crystal $\B(\infty)\otimes \B(1)\otimes\B(2)$
given by
\begin{equation*}
\mathcal{I}(\infty)=\{u_{\infty}\otimes\beta_1\otimes\beta_2\},
\end{equation*}
where
\begin{align}
\beta_1 &=b_1(-k_{1,\bar 2})\otimes
          b_2(-k_{1,\bar 3})\otimes
          b_1(-k_{1,3})\otimes
          b_2(-k_{1,2})\otimes
          b_1(-k_{1,1})\in \B(1),\\
\beta_2 &=b_2(-k_{2,2})\in \B(2),
\end{align}
and where $k_{u,v}$ are any nonnegative integers such that
\begin{equation*}
0\le k_{1,\bar 2}\le k_{1,\bar 3}\le k_{1,3}/2
 \le k_{1,2}\le k_{1,1}.
\end{equation*}
The set $\mathcal{I}(\infty)$ forms a $\uq(G_2)$-subcrystal of
$\B(\infty)\otimes \B(1)\otimes\B(2)$.
\end{prop}
\begin{proof}
It suffices to show that
the action of Kashiwara operators satisfy the following properties\;:
\begin{equation*}
\fit \mathcal{I}(\infty) \subset \mathcal{I}(\infty),\quad \eit
\mathcal{I}(\infty) \subset \mathcal{I}(\infty) \cup \{0\},
\end{equation*}
for all $i\in I$.

We will compute the value $\fit$ on each element of
$\mathcal{I}(\infty)$, using the tensor product rule given in
Proposition~\ref{propq}. First, we compute the finite sequence
$\{a_k\}$ set by~\eqref{kashi} for
\begin{align*}
b&=u_{\infty}\otimes\beta_1\otimes\beta_2\\
 &=u_{\infty}\otimes
  b_1(-k_{1,\bar 2})\otimes
  b_2(-k_{1,\bar 3})\otimes
  b_1(-k_{1,3})\otimes
  b_2(-k_{1,2})\otimes
  b_1(-k_{1,1})\otimes
  b_2(-k_{2,2}).
\end{align*}
In the $i=1$ case, we have
\begin{align*}
a_1&=0,\quad a_3=a_5=a_7=-\infty,\\
a_2&=k_{1,\bar 2},\quad a_4=k_{1,3}+2k_{1,\bar 2}-3k_{1,\bar 3},\\
a_6&=k_{1,1}+2k_{1,\bar 2}-3k_{1,\bar 3}+2k_{1,3}-3k_{1,2},
\end{align*}
and for $i=2$ case,
\begin{align*}
a_1&=0,\quad a_2=a_4=a_6=-\infty,\\
a_3&=k_{1,\bar 3}-k_{1,\bar 2},\quad a_5=k_{1,2}
   -k_{1,\bar 2}+2k_{1,\bar 3}-k_{1,3},\\
a_7&=k_{2,2}-k_{1,\bar 2}+2k_{1,\bar 3}-k_{1,3}+2k_{1,2}-k_{1,1}.
\end{align*}
By Proposition~\ref{propq}, we obtain the following three
candidates of $\tilde f_i(b)$ for each $i$\,:
\begin{align*}
\tilde f_1(b)&=u_{\infty}\otimes \tilde
   f_1(b_1(-k_{1,\bar 2}))\otimes b_2(-k_{1,\bar 3})
   \otimes b_1(-k_{1,3})\otimes b_2(-k_{1,2})\otimes
   b_1(-k_{1,1})\\
   &\qquad\ \ \otimes b_2(-k_{2,2})\\
&=u_{\infty}\otimes \big{(}b_1(-k_{1,\bar 2}-1)\otimes
b_2(-k_{1,\bar
   3})\otimes b_1(-k_{1,3})\otimes b_2(-k_{1,2})\otimes
   b_1(-k_{1,1})\big{)}\\
   &\qquad\ \ \otimes b_2(-k_{2,2})\in u_{\infty}\otimes
   \B(1)\otimes\B(2),
\end{align*}
when $a_2\ge a_k$ for $1\le k<2$ and $a_2>a_k$ for $2<k\le 7$,
\begin{align*}
\tilde f_1(b)&=u_{\infty}\otimes
   b_1(-k_{1,\bar 2})\otimes b_2(-k_{1,\bar 3})
   \otimes \tilde f_1(b_1(-k_{1,3})) \otimes
   b_2(-k_{1,2})\otimes b_1(-k_{1,1})\\
   &\qquad\ \ \otimes b_2(-k_{2,2})\\
&=u_{\infty}\otimes \big{(}b_1(-k_{1,\bar 2})\otimes
b_2(-k_{1,\bar 3})
   \otimes b_1(-k_{1,3}-1) \otimes b_2(-k_{1,2})\otimes
   b_1(-k_{1,1})\big{)}\\
   &\qquad\ \ \otimes b_2(-k_{2,2})\in u_{\infty}\otimes
   \B(1)\otimes\B(2),
\end{align*}
when $a_4\ge a_k$ for $1\le k<4$ and $a_4>a_k$ for $4<k\le 7$,
\begin{align*}
\tilde f_1(b)&=u_{\infty}\otimes
   b_1(-k_{1,\bar 2})\otimes b_2(-k_{1,\bar 3})
   \otimes b_1(-k_{1,3})\otimes b_2(-k_{1,2})\otimes
   \tilde f_1(b_1(-k_{1,1}))\\
   &\qquad\ \ \otimes b_2(-k_{2,2})\\
&=u_{\infty}\otimes \big{(}b_1(-k_{1,\bar 2})\otimes
   b_2(-k_{1,\bar 3})
   \otimes b_1(-k_{1,3}) \otimes b_2(-k_{1,2})\otimes
   b_1(-k_{1,1}-1)\big{)}\\
   &\qquad\ \ \otimes b_2(-k_{2,2})\in u_{\infty}\otimes
   \B(1)\otimes\B(2),
\end{align*}
when $a_6\ge a_k$ for $1\le k<6$ and $a_6>a_k$ for $6<k\le 7$,
\begin{align*}
\tilde f_2(b)&=u_{\infty}\otimes
   b_1(-k_{1,\bar 2})\otimes \tilde f_2(b_2(-k_{1,\bar 3}))
   \otimes b_1(-k_{1,3})\otimes b_2(-k_{1,2})\otimes
   b_1(-k_{1,1})\\
   &\qquad \ \ \otimes b_2(-k_{2,2})\\
&=u_{\infty}\otimes \big{(}b_1(-k_{1,\bar 2})\otimes
   b_2(-k_{1,\bar 3}-1)
   \otimes b_1(-k_{1,3})\otimes b_2(-k_{1,2})\otimes
   b_1(-k_{1,1})\big{)}\\
   &\qquad\ \ \otimes b_2(-k_{2,2})\in u_{\infty}\otimes
   \B(1)\otimes\B(2),
\end{align*}
when $a_3\ge a_k$ for $1\le k<3$ and $a_3>a_k$ for $3<k\le 7$,
\begin{align*}
\tilde f_2(b)&=u_{\infty}\otimes
   b_1(-k_{1,\bar 2})\otimes b_2(-k_{1,\bar 3})
   \otimes b_1(-k_{1,3}) \otimes \tilde
   f_2(b_2(-k_{1,2}))\otimes
   b_1(-k_{1,1})\\
   &\qquad\ \ \otimes b_2(-k_{2,2})\\
&=u_{\infty}\otimes \big{(}b_1(-k_{1,\bar 2})\otimes
   b_2(-k_{1,\bar 3})
   \otimes b_1(-k_{1,3}) \otimes b_2(-k_{1,2}-1)\otimes
   b_1(-k_{1,1})\big{)}\\
   &\qquad \ \ \otimes b_2(-k_{2,2})\in u_{\infty}\otimes
   \B(1)\otimes\B(2),
\end{align*}
when $a_5\ge a_k$ for $1\le k<5$ and $a_5>a_k$ for $5<k\le 7$,
\begin{align*}
\tilde f_2(b)&=u_{\infty}\otimes b_1(-k_{1,\bar
   2})\otimes b_2(-k_{1,\bar 3})
   \otimes b_1(-k_{1,3})\otimes b_2(-k_{1,2})\otimes
   b_1(-k_{1,1})\\
   &\qquad\ \ \otimes \tilde f_2(b_2(-k_{2,2}))\\
&=u_{\infty}\otimes \big{(}b_1(-k_{1,\bar 2})\otimes
   b_2(-k_{1,\bar 3})\otimes
   b_1(-k_{1,3}) \otimes b_2(-k_{1,2})\otimes
   b_1(-k_{1,1})\big{)}\\
   &\qquad\ \ \otimes b_2(-k_{2,2}-1)\in u_{\infty}\otimes
   \B(1)\otimes\B(2),
\end{align*}
when $a_7\ge a_k$ for $1\le k<7$. And for each case given above, we
obtain the following result from conditions for the sequence $a_k$. In
the $i=1$ case, $k_{u,v}$ values, appearing in the above expression for
$\fit b$, are nonnegative integers satisfying
\begin{itemize}
\item $k_{1,\bar 2}+1\le k_{1,\bar 3}\le k_{1,3}/2
      \le k_{1,2}\le k_{1,1}$,\\
      when $a_2\ge a_k$ for $1\le k<2$ and $a_2>a_k$ for $2<k\le 7$,
\item $k_{1,\bar 2}\le k_{1,\bar 3}\le (k_{1,3}+1)/2
      \le k_{1,2}\le k_{1,1}$,\\
      when $a_4\ge a_k$ for $1\le k<4$ and $a_4>a_k$ for $4<k\le 7$,
\item $k_{1,\bar 2}\le k_{1,\bar 3}\le k_{1,3}/2
      \le k_{1,2}\le k_{1,1}+1$,\\
      when $a_6\ge a_k$ for $1\le k<6$ and $a_6>a_k$ for $6<k\le 7$,
\end{itemize}
and in the $i=2$ case,
\begin{itemize}
\item $0\le k_{1,\bar 2}\le k_{1,\bar 3}+1\le k_{1,3}/2
      \le k_{1,2}\le k_{1,1}$,\\
      when $a_3\ge a_k$ for $1\le k<3$ and $a_3>a_k$ for $3<k\le 7$,
\item $0\le k_{1,\bar 2}\le k_{1,\bar 3}\le k_{1,3}/2
      \le k_{1,2}+1\le k_{1,1}$,\\
      when $a_5\ge a_k$ for $1\le k<5$ and $a_5>a_k$ for $5<k\le 7$,
\item $0\le k_{2,2}+1$,\\
      when $a_7\ge a_k$ for $1\le k<7$.
\end{itemize}
Thus the action of Kashiwara operator $\fit$ is closed on
$\mathcal{I}(\infty)$.

Proof for the statements concerning $\eit$
may be done in a similar manner.
\end{proof}
The notation $\beta_1$ and $\beta_2$ appearing in this proposition will
be used a few more times in this section.

\begin{theorem}\label{main2}
There exists a $\uq(G_2)$-crystal isomorphism
\begin{equation}
\B(\infty) \overset{\sim} \longrightarrow \T(\infty)
\overset{\sim} \longrightarrow \mathcal{I}(\infty)
\subset\B(\infty)\otimes \B(1)\otimes\B(2),
\end{equation}
which maps $T_{\infty}$ to $u_\infty\otimes
 (b_1\otimes b_2\otimes b_1\otimes b_2\otimes b_1)\otimes (b_2)$.
\end{theorem}
\begin{proof}
With the help of tensor product rules, it is easy to check the
compatibility of this map with Kashiwara operators. Other parts of
the proof are similar or easy. Hence we shall only write out the
maps and give no proofs.

For each tableau with the second row consisting of $b^2_3$-many
$3$-boxes and just one $2$-box, and with the first row consisting of
$b^1_j$-many $j$-boxes, for each $j\succ 1$, and $(b^2_3+2)$-many
$1$-boxes, we may map it to the element
$u_{\infty}\otimes\beta_1\otimes\beta_2$ where
\begin{align*}
k_{1,1}&=\sum_{j=2}^{\bar 1} b^1_j,\quad k_{1,2}=\sum_{j=3}^{\bar
1} b^1_j,\quad k_{1,3}=2(\sum_{j=\bar 3}^{\bar 1} b^1_j)+b_0^1,\\
k_{1,\bar 3}&=b^1_{\bar{2}}+b^1_{\bar 1},\quad k_{1,\bar
2}=b^1_{\bar 1},\quad k_{2,2}=b^2_3.
\end{align*}

Conversely, an element $u_{\infty}\otimes\beta_1\otimes\beta_2$ is
sent to the tableau whose shape we describe below row-by-row.
\begin{itemize}\allowdisplaybreaks
\item The first row consists of
  \begin{alignat*}{2}
  &\text{$(k_{1,\bar{2}})$-many $\bar 1$\,s,}&
  &\text{$(k_{1,\bar{3}}-k_{1,\bar{2}})$-many $\bar 2$\,s,}\\
  &\text{$\lfloor k_{1,3}/2 - k_{1,\bar 3}\rfloor$-many $\bar 3$\,s,}&
  &\text{($(A+B)-(A'+B')$)-many $0$\,s,}\\
  &\text{$(k_{1,2}-k_{1,3}/2)\)-many $3$\,s,}&\
  &\text{($k_{1,1}-k_{1,2}$)-many $2$\,s, and}\\
  &\text{$\textstyle
  \big(k_{2,2}+2\big)$-many $1$\,s.}
  \end{alignat*}
\item The second row consists of
      \begin{align*}
      \text{$(k_{2,2})$-many $3$\,s\quad and\quad one $2$}.
      \end{align*}
\end{itemize}
Here, $A=k_{1,2}-k_{1,3}/2$, $B=k_{1,3}/2-k_{1,\bar 3}$, $A'=\lfloor
k_{1,2}-k_{1,3}/2]$, and $B'=\lfloor k_{1,3}/2-k_{1,\bar 3}\rfloor$.
\end{proof}

Since the above theorem has shown $\B(\infty)\cong \mathcal{I}(\infty)$
as crystals, image of the injective crystal morphism
\begin{equation*}
\Psi : \Binf\rightarrow\Binf\otimes \B(1)\otimes \B(2)= \Binf\otimes
   (\B_1\otimes\B_2\otimes\B_1\otimes\B_2\otimes\B_1)\otimes(\B_2),
\end{equation*}
which maps $u_{\infty}$ to $u_\infty\otimes
 (b_1\otimes b_2\otimes b_1\otimes b_2\otimes b_1)\otimes (b_2)$
is $\mathcal{I}(\infty)$.

In the following corollary, a combinatorial description of
$\B(\infty)$ for $G_2$-type is given following Cliff's method. A
specific choice for the index sequence of crystals $S =
(1,2,1,2,1,2)$ corresponding to a longest word
$w_0=s_{1}s_{2}s_{1}s_{2}s_{1}s_{2}$ of the Weyl group is used.

\begin{cor}
Image of the injective strict crystal morphism
\begin{equation*}
\Psi : \Binf\rightarrow\Binf\otimes
   (\B_1\otimes\B_2\otimes\B_1\otimes\B_2\otimes\B_1)\otimes(\B_2),
\end{equation*}
which maps $u_{\infty}$ to $u_\infty\otimes
 (b_1\otimes b_2\otimes b_1\otimes b_2\otimes b_1)\otimes (b_2)$
is given by
\begin{equation*}
\Psi(\Binf) = \mathcal{I}(\infty)=\{ u_\infty \otimes \beta_1
\otimes \beta_2 \}.
\end{equation*}
\end{cor}

We illustrate the correspondence between $\T(\infty)$ and
$\Psi(\B(\infty))$ for type $G_2$.

\begin{example}
The marginally large tableau
\begin{center}
\begin{texdraw}%
\drawdim in \arrowheadsize l:0.065 w:0.03 \arrowheadtype t:F
\fontsize{5}{5}\selectfont \textref h:C v:C \drawdim em
\setunitscale 1.6 \move(0 2)\rlvec(9 0) \move(0 1)\rlvec(9
0)\rlvec(0 1) \move(0 0)\rlvec(3 0)\rlvec(0 2) \move(0 2)\rlvec(0
-2) \move(1 0)\rlvec(0 2) \move(2 0)\rlvec(0 2) \move(4 1)\rlvec(0
1) \move(5 1)\rlvec(0 1) \move(6 1)\rlvec(0 1) \move(7 1)\rlvec(0
1) \move(8 1)\rlvec(0 1) \htext(0.5 1.5){$1$} \htext(1.5 1.5){$1$}
\htext(2.5 1.5){$1$} \htext(3.5 1.5){$1$} \htext(4.5 1.5){$2$}
\htext(5.5 1.5){$0$} \htext(6.5 1.5){$\bar 3$} \htext(7.5
1.5){$\bar 3$} \htext(8.5 1.5){$\bar 1$} \htext(0.5 0.5){$2$}
\htext(1.5 0.5){$3$} \htext(2.5 0.5){$3$}
\end{texdraw}%
\end{center}
of $\T(\infty)$ corresponds to the element
\begin{equation*}
u_{\infty}\otimes
  b_1(-1)\otimes
  b_2(-1)\otimes
  b_1(-7)\otimes
  b_2(-4)\otimes
  b_1(-5)\otimes
  b_2(-2)
\end{equation*}
of $\Psi(\B(\infty))$ under the map given in Theorem~\ref{main2}.
\end{example}

\begin{remark}
We gave two new explicit descriptions of the crystal $\Binf$ in
this paper, namely, $\Psi(\Binf)$ and $\M(\infty)$ for the $G_2$
case. We can provide maps between the two giving crystal
isomorphisms
\begin{equation*}
 \M(\infty)\overset {\sim} \longrightarrow\Psi(\B(\infty))
\end{equation*}
in both directions. The maps can easily be drawn from
Theorem~\ref{mainreal5} and Theorem~\ref{main2}.
\end{remark}

%%%%%%%%%%%%%%%%%%%%%%%%%%%%%%%%%%%%%%%%%%%%%%%%%%%%%%%%%%%%%%%%%%%%%%%%%%
\bibliographystyle{amsplain}

\begin{thebibliography}{10}

\bibitem{MR1614241}
G. Cliff, \emph{Crystal bases and {Y}oung tableaux}, J. Algebra
  \textbf{202} (1998), no.~1, 10--35.

\bibitem{MR1881971}
J. Hong and S.-J. Kang, \emph{Introduction to quantum groups and crystal
  bases}, Graduate Studies in Mathematics, vol.~42, Amer. Math.
  Soc., Providence, RI, 2002.

\bibitem{HL}
J. Hong and H. Lee, \emph{Young tableaux and crystal $\Binf$ for
finite simple Lie algebras}, arXiv:math.QA/0507448.

\bibitem{MR1359532}
J.~C. Jantzen, \emph{Lectures on quantum groups}, Graduate Studies in
  Mathematics, vol.~6, Amer. Math. Soc., Providence, RI, 1996.

\bibitem{MR2043368}
S.-J. Kang, J.-A. Kim, and D.-U. Shin,
  \emph{Monomial realization of crystal bases for special linear {L}ie
  algebras}, J. Algebra \textbf{274} (2004), no.~2, 629--642.

\bibitem{MR2074700}
\bysame, \emph{Crystal bases for quantum
  classical algebras and {N}akajima's monomials}, Publ. Res. Inst. Math. Sci.
  \textbf{40} (2004), no.~3, 757--791.

\bibitem{KM}
S.-J. Kang and K. C. Misra, \emph{Crystal bases and tensor product
decompositions of {$U\sb q(G\sb 2)$}-modules}, J. Algebra
\textbf{163} (1994), no.~3, 675--691.

\bibitem{K90}
M.~Kashiwara, \emph{Crystalizing the $q$-analogue of universal
enveloping algebras} Comm. Math. Phys. \textbf{133} (1990), no.~2,
249--260.

\bibitem{K91}
\bysame, \emph{On crystal bases of the {$Q$}-analogue of universal
enveloping algebras}, Duke Math. J. \textbf{63} (1991), no.~2,
465--516.

\bibitem{MR1240605}
\bysame, \emph{The crystal base and {L}ittelmann's refined
{D}emazure character formula}, Duke Math. J. \textbf{71} (1993),
no.~3, 839--858.

\bibitem{MR1988989}
\bysame, \emph{Realizations of crystals}, Combinatorial and
geometric
  representation theory (Seoul, 2001), Contemp. Math., vol. 325,
  Amer. Math. Soc., Providence, RI, 2003, pp.~133--139.

\bibitem{ks}
M. Kashiwara and Y. Saito, \emph{Crystal graphs for
  representations of the {$q$}-analogue of classical {L}ie algebras}, J.
  Algebra \textbf{165} (1994), no.~2, 295--345.

\bibitem{kim}
J.-A. Kim, \emph{Monomial realization of crystal graphs for
  $\uq(A_n^{(1)})$}, to appear in Math. Annel.

\bibitem{lee2}
H. Lee, \emph{Young tableaux, Nakajima monomials, and crystals for
special linear Lie algebras}, arXiv:math.QA/0506147.

\bibitem{lee}
\bysame, \emph{Extended Nakajima's monomials and realizations of
crystals for $\hat{\mathfrak{sl}}_n$}, KIAS preprint M05008.

\bibitem{MR1865400}
H. Nakajima, \emph{Quiver varieties and tensor products}, Invent. Math.
  \textbf{146} (2001), no.~2, 399--449.

\bibitem{MR1988990}
\bysame, \emph{{$t$}-analogs of {$q$}-characters of quantum affine algebras of
  type {$A\sb n,D\sb n$}}, Combinatorial and geometric representation theory
  (Seoul, 2001), Contemp. Math., vol. 325, Amer. Math. Soc., Providence, RI,
  2003, pp.~141--160.

\bibitem{MR1475048}
T. Nakashima and A. Zelevinsky, \emph{Polyhedral realizations of
  crystal bases for quantized {K}ac-{M}oody algebras}, Adv. Math. \textbf{131}
  (1997), no.~1, 253--278.

\bibitem{saito}
Y. Saito, \emph{Combinatorial and geometric realization of
  crystal ${\B}(\infty)$ for type ${A}_n$}, KIAS lecture series 2002.

\bibitem{shin}
D.-U. Shin,
  \emph{Crystal bases and monomials for $U_q(G_2)$-modules},
  arXiv:math.QA/0308176 .

\end{thebibliography}

\providecommand{\bysame}{\leavevmode\hbox to3em{\hrulefill}\thinspace}
\providecommand{\MR}{\relax\ifhmode\unskip\space\fi MR }
\providecommand{\MRhref}[2]{%
  \href{http://www.ams.org/mathscinet-getitem?mr=#1}{#2}
}
\providecommand{\href}[2]{#2}

%%%%%%%%%%%%%%%%%%%%%%%%%%%%%%%%%%%%%%%%%%%%%%%%%%%%%%%%%%%%%%%%%%%%%%%%%%
\end{document}